\documentclass[12pt,reqno]{article}
\usepackage{amsmath}
\usepackage{amssymb}
\usepackage{amsthm}
\usepackage{mathrsfs}
\usepackage{latexsym}
\usepackage{color}
\usepackage{authblk}

\newtheorem{theorem}{Theorem}[section]
\newtheorem{lemma}{Lemma}[section]
\newtheorem{corollary}{Corollary}[section]
\newtheorem{remark}{Remark}[section]
\newtheorem{definition}{Definition}[section]
\newtheorem{proposition}{Proposition}[section]
\newtheorem{example}{Example}[section]
\newtheorem{assumption}{Assumption}[section]
\numberwithin{equation}{section}
\newcommand{\bth}{\begin{theorem}}
\newcommand{\ethe}{\end{theorem}}

\newcommand{\bre}{\begin{remark}}
\newcommand{\ere}{\end{remark}}

\newcommand{\ble}{\begin{lemma}}
\newcommand{\ele}{\end{lemma}}
\newcommand{\bde}{\begin{definition}}
\newcommand{\ede}{\end{definition}}

\newcommand{\bco}{\begin{corollary}}
\newcommand{\eco}{\end{corollary}}

\newcommand{\bpr}{\begin{proposition}}
\newcommand{\epr}{\end{proposition}}

\newcommand{\bexer}{\begin{exercise}}
\newcommand{\eexer}{\end{exercise}}
\newcommand{\breh}{\begin{hint}}
\newcommand{\ereh}{\end{hint}}

\newcommand{\halmos}{\hfill \qed}

\newcommand{\bexam}{\begin{example}}
\newcommand{\eexam}{\end{example}}

\newcommand{\pr} {{\bf Proof.}}

\newcommand{\bfi}{\begin{fig}}
\newcommand{\efi}{\end{fig}}

\newcommand{\beao}{\begin{eqnarray*}}
\newcommand{\eeao}{\end{eqnarray*}\noindent}

\newcommand{\beam}{\begin{eqnarray}}
\newcommand{\eeam}{\end{eqnarray}\noindent}
\newcommand{\E}{\mathbf{E}}
\newcommand{\PP}{\mathbf{P}}

\newcommand{\xto}{x\to\infty}

\newcommand{\bF}{\overline{F}}

\newcommand{\bG}{\overline{G}}
\newcommand{\bV}{\overline{V}}

\newcommand{\bbr}{{\mathbb R}}

\newcommand{\bbn}{{\mathbb N}}

\newcommand{\yto}{y\to\infty}

\newcommand{\vep}{\varepsilon}

\begin{document}
\title{Asymptotics for aggregated interdependent multivariate subexponential claims with general investment returns}

\author[1]{\small Zhangting Chen\thanks{\texttt{20234007012@stu.suda.edu.cn}}}

\author[2]{Dimitrios G. Konstantinides\thanks{\texttt{konstant@aegean.gr}}}

\author[2]{Charalampos  D. Passalidis\thanks{\texttt{sasd24009@sas.aegean.gr}}}

\affil[1]{\small School of Mathematical Sciences, Soochow University, Suzhou, Jiangsu, China}
\affil[2]{Dept. of Statistics and Actuarial-Financial Mathematics, University of the Aegean, Karlovassi, GR-83 200 Samos, Greece}

\date{}
\maketitle
\begin{abstract}
This paper investigates asymptotic estimates for the entrance probability of the discounted aggregate claim vector from a multivariate renewal risk model into some rare set. We provide asymptotic results for the entrance probability on both finite and infinite time horizons under various assumptions regarding the stochastic price process of the investment portfolio, the distribution class of claim vectors, and the dependence structure among the claim vectors. We note that the main results extend beyond the class of multivariate regular variation. Furthermore, we introduce two dependence structures to model the dependence among the claim vectors. In particular, our results are new even in the one-dimensional subcase.
\end{abstract}

\maketitle
\vspace{3mm}



\section{Introduction} \label{sec.CCKP.1}

Consider an insurance company that operates $d$-lines of business, with $d \in \bbn$, and simultaneously makes risk-free and risky investments. Consequently, the insurer is exposed to two fundamental risks, which are insurance risk due to insurance claims and financial risk caused by investments. We may find recently a plethora of papers on risk theory which study risk models with insurance and financial risks, either on discrete time or on continuous time, see for example \cite{tang:tsitsiashvili:2003a}, \cite{li:tang:2015}, \cite{yang:konstantinides:2015}, \cite{chen:yuan:2017}, \cite{tang:yang:2019}, \cite{cui:wang:2024} and \cite{paulsen:gjessing:1997}, \cite{paulsen:2002}, \cite{tang:wang:yuen:2010}, \cite{li:2012}, \cite{guo:wang:2013}, \cite{yang:wang:konstantinides:2014}, respectively. This paper focuses on the continuous-time case. In this model, the insurer is faced with successive claim vectors $\{{\bf X}^{(i)}\,,\;i \in \bbn\}$ that happened at the random times $\{\tau_i\,,\;i \in \bbn\}$. And the random times form a renewal counting process%
\begin{align*}
N(t) :=\sup\{i\in \bbn\;:\;\tau_i \leq t\}\,,\,\, t\geq 0\,,
\end{align*}%
with convention $\sup \emptyset =0$, and finite renewal function $\lambda(t):=\E[N(t)]=\sum_{i=1}^{\infty}\PP[\tau_i \leq t]$. We note that, throughout the paper, the claim vectors  $\{{\bf X}^{(i)}\,,\; i \in \bbn\}$ are identically distributed copies of ${\bf X}$ with distribution $F$, whose support is restricted to the nonnegative quadrant, and each vector ${\bf X}^{(i)}=(X_1^{(i)},\,\ldots,\, X_d^{(i)})^{\top}$ (where ${\bf x}^{\top}$ denotes the transpose of vector ${\bf x}$) may contain zero components but cannot be a zero vector. For papers on bivariate or multivariate risk models, we refer the reader to \cite{chen:li:cheng:2023}, \cite{li:2023}, \cite{yang:chen:yuen:2024}, \cite{xu:shen:wang:2025}, \cite{li:2016}, \cite{cheng:konstantinides:wang:2022}, \cite{yang:su:2023}, \cite{konstantinides:passalidis:2024j}.

We consider that $F$ is a multivariate heavy-tailed distribution, such as multivariate subexponential or multivariate consistently varying (see definitions in Section 2 below). These multivariate distribution classes are more general than the well-known multivariate regularly varying (${MRV}$) distribution class. For results on multivariate risk models with ${MRV}$ claims, the reader is referred to \cite{konstantinides:li:2016}, \cite{li:2016}, \cite{cheng:konstantinides:wang:2022}, and \cite{yang:su:2023}, and further, those papers assumed that each claim vector has asymptotically dependent components.  In our framework, we also allow many asymptotically independent structures. Moreover, this paper allows certain dependence structure, that incude independence as a special case, among the claim vectors and thus extends beyond the independence assumption of $\{{\bf X}^{(i)}\,,\; i \in \bbn\}$ used in prior literature. This approach establishes interdependence among the claim vectors, namely, each vector has dependent components, while simultaneously the vectors are dependent on each other. To the best of our knowledge, the continuous-time risk models with interdependent claims have been studied only by \cite{chen:cheng:zheng:2025}, \cite{konstantinides:passalidis:2024j}. More discussion about interdependence can be found in \cite{konstantinides:passalidis:2023} and \cite{konstantinides:passalidis:2025}.

We suppose that the insurer invests its surplus in risk-free and risky assets, and the logarithmic returns of the investment portfolio are depicted by a general c\'{a}dl\'{a}g process $\{\xi(t)\,,\;t\geq 0\}$ with $\xi(0)=0$. Then, for any fixed $T \in \left\{t\geq 0:0<\lambda (t)\leq \infty\right\}=:\varLambda$, the insurer's discounted aggregated claim vector is given through the relation
\beam \label{eq.CCKP.1.2}
{\bf D}(T) =\sum_{i=1}^{N(T)} {\bf X}^{(i)}\,e^{-\xi(\tau_i)}=\left( %
\begin{array}{c}
\sum_{i=1}^{N(T)} X_{1}^{(i)}\,e^{-\xi(\tau_{i}) } \\ 
\vdots \\ 
\sum_{i=1}^{N(T)} X_{d}^{(i)}\,e^{-\xi(\tau_{i}) } 
\end{array} 
\right)\,,
\eeam
where $N(\infty)$ is understood $\infty$. Next, we present the basic assumptions for the whole paper.

\begin{assumption} \label{ass.CCKP.1.1}
The vectors $\{{\bf X}^{(i)}\,,\;i \in \bbn\}$, are identically distributed with common distribution $F$, which has support on the nonnegative quadrant. Further, $\{{\bf X}^{(i)}\,,\;i \in \bbn\}$, $\{N(t)\,,\;t\geq 0 \}$ and  $\{\xi(t)\,,\;t \geq 0\}$, are mutually independent. 
\end{assumption} 

\begin{assumption}\label{ass.CCKP.1.2}
	$\{\xi(t)\,,\;t \geq 0\}$ is a general stochastic process bounded away from $-\infty$, i.e., for any $T\in\varLambda\setminus\{\infty\}$, there exists some positive constant $K_T$ such that 
\begin{align}\label{th.CZT.3.1}
\PP\left[\inf_{0\leq t\leq T} \xi(t) \geq -K_T\right]=1\,.
\end{align}
\end{assumption}
Relation (\ref{th.CZT.3.1}) means that the process $\{e^{\xi(t)}\, , \; 0\leq t\leq T\}$ is away from zero. As noted by \cite{chen:liu:2023}, the admissibility of an investment strategy for a self-financing portfolio in portfolio theory is often defined by conditions (\ref{th.CZT.3.1}), which ensure realistic and sustainable trading behavior (see, e.g., \cite[Def. 10.2]{Bjork:2009}). In the general financial markets, investors can take short positions, allowing asset prices to theoretically become negative. However, the insurance sector operates under stricter regulatory frameworks, such as those imposed by the National Association of Insurance Commissioners, which largely prohibit insurers from engaging in short selling. For further discussion on this topic, see \cite{Molk:Partnoy:2019}, who analyze the role of institutional investors, particularly insurers, in short selling. Their study reveals that insurers rarely engage in short selling, and when they do, it is primarily for hedging rather than speculative purposes. This aligns with condition (\ref{th.CZT.3.1}) in insurance settings, where restrictive regulatory frameworks discourage aggressive short positions. The reader is referred to \cite[Ex. 1]{chen:liu:2023}, which shows that virtually all investment strategies in insurance practice should fulfill condition (\ref{th.CZT.3.1}).

This paper investigates the asymptotic behavior of the probability $\PP\left[ {\bf D}(T) \in x\, A  \right]$ as $x\to\infty$, where $x\, A$ is some rare set, coming from a rather general set family capable of representing insurance-relevant scenarios. Namely, we establish the following asymptotic relation
\beam \label{eq.CCKP.1.3}
\lim_{x\to\infty}\frac{\PP\left[ {\bf D}(T) \in x\,A  \right]}{\int_0^T \PP[{\bf X}\,e^{-\xi(s)} \in x\,A] \,\lambda(ds)}\,=\,1\,,
\eeam
for any fixed $T \in\varLambda$. As mentioned previously, our asymptotic results drop the condition of the independence assumption of $\{{\bf X}^{(i)}\,,\; i \in \bbn\}$ used in prior literature and allow some weak dependence among the claim vectors. The logarithmic return of the investment portfolio $\{\xi(t)\,,\; t\geq 0\}$ is a general process satisfying some mild conditions that are satisfied by many stochastic processes in financial mathematics. 

The rest of this paper is organized as follows. In Section 3, we give three theorems, together with some corollaries, that establish the relation \eqref{eq.CCKP.1.3} with $T<\infty$ or $T=\infty$ under the assumptions that distribution $F$ belongs to multivariate subexponential class, the vectors $\{{\bf X}^{(i)}\,,\;i \in \bbn\}$ satisfy a weak dependence condition, and the $\{\xi(t)\,,\;t\geq 0\}$ is a general stochastic process.
In Section 4, we provide the proofs of the main results, together with the corresponding auxiliary lemmas.

\section{Preliminaries} \label{sec.CCKP.2}

Now, we recall some definitions of multivariate heavy-tailed distributions, as well as some dependence structures. Before this, some notations are provided.

Throughout the paper, the asymptotic expressions are understood in the sense of $\xto$, except otherwise stated. All the vectors are depicted with bold script and are of dimension $d \in \bbn$. For any two real numbers $a$ and $b$, we write $\max\{a,\,b\}=a\vee b$ and $\min\{a,\,b\}=a\wedge b$ and with $\left\lfloor a \right\rfloor$ we denote the integer part of $a$. Furthermore, by convention, the $\sum_{i\in E}=0$, when $E=\emptyset$. For any two vectors ${\bf x},\,{\bf y} $, we denote ${\bf x}\pm{\bf y} =(x_1\pm y_1,\,\ldots,\,x_d\pm y_d)^{\top}$, and the scalar product, with some positive quantity $c>0$, is represented by $c\,{\bf x} =(c\,x_1,\,\ldots,\,c\,x_d)^{\top}$. Furthermore, we use a bold Arabic number to denote the vector with all elements being that number, e.g., ${\bf 0}=(0,\,\ldots,\,0)^{\top}$ for the origin of the axes. For any set $B$ from space $\bbr^d$, we denote by $B^c$ its complement set, by $\partial B$ its border, by $\overline{B}$ its closed hull, and by ${\bf 1}_{B}$ its indicator function. We say that the set  $B$ is ``increasing" if for any ${\bf x} \in B$ and any  ${\bf y} \in \bbr_{+}^d :=[0,\,\infty)^d$, it holds ${\bf x} +{\bf y} \in B$. 

For two positive functions ${f},\,{g}$, we write ${f}(x) \sim {g}(x)$, if $\lim \frac{{f}(x)}{{g}(x)}=1$; write ${f}(x) =o[ {g}(x)]$, if $\lim \frac{{f}(x)}{{g}(x)}= 0$; write ${f}(x) \lesssim {g}(x)$ (or, equivalently, ${g}(x) \gtrsim {f}(x) $), if $\limsup \frac{{f}(x)}{{g}(x)}\leq 1$; write  ${f}(x) =O[ {g}(x)]$, if $\limsup \frac{{f}(x)}{{g}(x)}< \infty$. Further, we write ${ f}(x) \asymp {g}(x)$ if both ${f}(x)=O[{g}(x)]$ and ${g}(x)=O[{f}(x)]$.

For any distribution $G$, we denote the tail of $G$ by $\bG(x)=1-G(x)$, for any $x \in \bbr$. For two distributions $G_1$, $G_2$ we denote by $G_1\ast G_2$, 
their convolution while we denote by $G^{n\ast}$ the $n$-th convolution power of $G$ with itself.

\subsection{Heavy-tailed distributions}

To introduce the multivariate distribution classes, we need at first the one-dimensional ones. In what follows in this subsection, we assume that for a univariate distribution $G$ supported on $\bbr_+$, it holds that $\bG(x)>0$ for any $x \ge 0$.

One of the largest heavy-tailed distribution classes is the class $\mathcal{L}$, of long-tailed distributions. It is said that $G \in \mathcal{L}$, if for any (or, equivalently, for some) $y>0$ it holds that
\begin{align*}
{\bG(x-y)}\sim{\bG(x)}\,.
\end{align*}
One of the most famous subclasses of $\mathcal{L}$ is the subexponential distribution class, symbolically $\mathcal{S}$. It is said that $ G\in\mathcal{S}$, if for any (or, equivalently, for some) integer $n\geq 2$ it holds that
\begin{align*}
{\overline{G^{n\ast}}(x)}\sim n{\bG(x)}\,.
\end{align*}
Classes $\mathcal{L}$ and $\mathcal{S}$ were introduced in \cite{chistyakov:1964}, where we find also the inclusion $\mathcal{S} \subsetneq \mathcal{L}$. Class $\mathcal{S}$ finds many applications in several branches of applied probability, as for example in risk theory, queuing theory, random walks, among others, see in \cite{embrechts:klueppelberg:mikosch:1997}, \cite{asmussen:2003}, \cite{borovkov:borovkov:2008} for related monographs.

Another important heavy-tailed distribution class is the class $\mathcal{D}$ of dominatedly varying distributions introduced in \cite{feller:1969}. It is said that $G \in \mathcal{D}$, if for any (or, equivalently, for some) $b \in (0,1)$ it holds that
\begin{align*}
\limsup \dfrac{\bG(b\,x)}{\bG(x)} < \infty\,.
\end{align*}
 It is known that $\mathcal{D} \not\subsetneq \mathcal{S}$ and $\mathcal{S} \not\subsetneq \mathcal{D}$ and further it holds $\mathcal{D}\cap \mathcal{S} \equiv \mathcal{D}\cap \mathcal{L}$, see in \cite{goldie:1978}. A smaller class than $\mathcal{D}\cap \mathcal{L}$, is the class of consistently varying distributions $\mathcal{C}$, introduced in \cite{cline:1994}. We say that $G\in \mathcal{C}$, if it holds that
\begin{align*}
\lim_{b \uparrow 1} \limsup \dfrac{\bG(b\,x)}{\bG(x)} =1\,.
\end{align*} 
Finally, we introduce a smaller distribution class called the regularly varying distribution class $\mathcal{R}_{-\alpha}$, with $\alpha \in (0,\,\infty)$. It is said that $G\in \mathcal{R}_{-\alpha}$, with index of regular variation $\alpha \in (0,\,\infty)$, if for any $y>0$, it holds
\begin{align*}
{\bG(y\,x)}\sim y^{-\alpha}{\bG(x)}\,.
\end{align*}  
Let us note that the following inclusions are true, see for example in \cite{foss:korshunov:zachary:2013}, or in \cite{leipus:siaulys:konstantinides:2023}:
\beam \label{eq.CCKP.2.6}
\mathcal{R}:=\bigcup_{0<\alpha < \infty}\mathcal{R}_{-\alpha} \subsetneq \mathcal{C} \subsetneq \mathcal{D}\cap \mathcal{L} \subsetneq \mathcal{S} \subsetneq \mathcal{L} \,.
\eeam 

The upper and lower Matuszewska indices are two well-known indices for characterizing heavy-tailed distributions. By definition, the upper and lower Matuszewska indices of a distribution $G$ are formulated as 
\begin{align*}
J_G^+:=-\lim_{\yto} \dfrac{\log \bG_*(y)}{\log y}\,,\qquad J_G^-:=-\lim_{\yto} \dfrac{\log \bG^*(y)}{\log y}\,,%
\end{align*}
respectively, where $\bG_*(y):=\liminf \frac{ \bG(y\,x)}{\bG(x)}$ and $\bG^*(y):=\limsup \frac{ \bG(y\,x)}{\bG(x)}\,$, for $y>1$. It was shown that $G \in \mathcal{D}$ if and only if $J_G^+<\infty$, and if $G \in \mathcal{R}_{-\alpha} $, with $\alpha \in (0,\,\infty)$, then $J_G^+ = J_G^- = \alpha$. Further, if $G \in\mathcal{D}$, then for any $p> J_G^+$, it holds that
\beam \label{eq.CCKP.2.8}
x^{-p}=o\left[\bG(x)\right]\,.
\eeam
For further discussions, we refer to \cite[Sec. 2.1.2]{bingham:goldie:teugels:1987} or \cite[Sec. 2.4]{leipus:siaulys:konstantinides:2023}. 

The following distribution class is immediately related to the lower Matuszeska index, and was introduced in \cite{haan:resnick:1984}. A distribution $G$ belongs to the class of positively decreasing distributions, symbolically $G \in \mathcal{P_D}$, if for any (or, equivalently, for some) $y>1$, it holds that
\begin{align*}
\bG^*(y)  < 1\,.
\end{align*}
It is known that $G \in \mathcal{P_D}$ if and only if $J_G^->0$. Furthermore, the class $ \mathcal{P_D}$ is wide enough to contain both light-tailed and heavy-tailed distributions. As mentioned by \cite[Sec. 2]{tang:2006}, almost all commonly used subexponential distributions are also in $\mathcal{P_D}$. Further, by the characterization via Matuszewska indices, it is easily implied that $\mathcal{R} \subsetneq \mathcal{C}\cap \mathcal{P_D}$, in combination with relation \eqref{eq.CCKP.2.6}. For further discussions about $\mathcal{P_D}$, we refer the reader to \cite{bardoutsos:konstantinides:2011} and \cite{konstantinides:passalidis:2025b}.

After establishing the concepts of one-dimensional heavy-tailed distributions, we consider the multivariate heavy-tailed distribution classes. Let us start with the ${MRV}$ distribution class, introduced in \cite{dehaan:resnick:1981}. 

Let ${\bf Z}$ be a random vector with distribution $V$, whose support is defined on the nonnegative quadrant. It is said that $V$ is ${MRV}$, denoted by $V \in {MRV}_{-\alpha}(\mu,\,G)$, if there exists some one-dimensional distribution $G \in \mathcal{R}_{-\alpha} $, with $\alpha \in (0,\,\infty)$, and some Radon measure $\mu$, nondegenerate to zero, such that
\begin{align}\label{ASSMRV}
\lim \dfrac{\PP[{\bf Z} \in x\,B]}{\bG(x)}=\mu(B)\,,
\end{align}
for every $\mu$-continuous Borel set $B \subsetneq [0,\,\infty]^d$ that is bounded away from ${\bf 0}$. Throughout this paper, when referring to the common distribution $F$ of claim vectors as ${MRV}$ without risk of confusion, we always assume that \eqref{ASSMRV} holds for ${\bf X}$, with ${\bf X}$ in place of ${\bf Z}$.
The ${MRV}$ distribution class has found many applications in time series, risk theory, and risk management, for example, we refer to \cite{resnick:2007}, \cite{konstantinides:li:2016}, \cite{samorodnitsky:2016}, \cite{li:2022b}, \cite{chen:liu:2024}, \cite{cheng:konstantinides:wang:2024}, \cite{liu:shushi:2024}, \cite{mikosch:wintenberger:2024}, \cite{chen:cheng:zheng:2025}, among others.

Although the ${MRV}$ distribution class is well-established, it still seems restrictive under some practical conditions. In the insurance scenario, modelling insurance claims by extremely heavy-tailed distributions may cause a significant overestimate of solvency capital requirements. In practice, moderately heavy-tailed claims may be more relevant. In \cite{samorodnitsky:sun:2016}, the multivariate subexponentiality was introduced by using the following family of sets.
\beam \label{eq.CCKP.2.10}
\mathscr{R}:=\left\{ A \subset \bbr^d\;:\; A\;{\text open,\, increasing},\;A^c\;{\text convex},\;{\bf 0}\notin \overline{A} \right\}\,.
\eeam
Some characteristic examples of sets, belonging to the family $\mathscr{R}$, that are interesting for the actuarial practice, are given as follows:
\beam \label{eq.CCKP.2.11}
A_1:=\left\{ {\bf y}\;:\; y_i>b_i\,,\;\exists\;i=1,\,\ldots,\,d \right\}\,,
\eeam
for $b_1,\,\ldots,\,b_d > 0$, and 
\beam \label{eq.CCKP.2.12}
A_2:=\left\{ {\bf y}\;:\; \sum_{i=1}^d l_i\,y_i>b \right\}\,,
\eeam
where $b>0$ and $l_1,\,\ldots,\,l_d \geq 0$, with 
\beao
\sum_{i=1}^d l_i = 1\,.
\eeao
Particularly, in the case of $d=1$, the set  $A_1$ is reduced to:
\beam \label{eq.CCKP.2.13}
A_3:=(b,\,\infty)\,,
\eeam
where $b>0$.

Hence, choosing some set from the \eqref{eq.CCKP.2.11}--\eqref{eq.CCKP.2.13}, the asymptotic estimation \eqref{eq.CCKP.1.3} acquires direct sense in actuarial practice. Indeed, if we choose $A_1$, then by \eqref{eq.CCKP.1.3} we get the asymptotic behavior of the probability, in one of the $d$ lines of business, the aggregated claims to exceed the initial capital that corresponds to it (or, the corresponding initial capital multiplied by some constant $b_i$, for $b_i \neq 1$). If we choose $A_2$, then by \eqref{eq.CCKP.1.3} we get the asymptotic behavior of the probability that the sum of aggregated claims of the $d$ lines of business exceeds the initial capital (or, the initial capital multiplied by some constant $b$, for $b \neq 1$). Finally, in case $d=1$, if we choose $A_3$, or even for $b=1$ in \eqref{eq.CCKP.2.13}, then by \eqref{eq.CCKP.1.3} we get the asymptotic behavior of the tail of discounted aggregated claims.

According to \cite[Lem. 4.5]{samorodnitsky:sun:2016}, the random variable
\beam \label{eq.CCKP.2.14}
Z_A:=\sup\left\{u\;:\;{\bf Z} \in u\,A \right\}\,,
\eeam
has a proper distribution $V_A$ given by
\beam \label{eq.CCKP.2.15}
\bV_A(x)=\PP[Z_A > x]=\PP[{\bf Z} \in x\,A]=\PP\left[\sup_{{\bf p}\in I_A}\;{\bf p}^{\top}\,{\bf Z} >x \right]\,,
\eeam
where $x>0$, see in \cite[Lem. 4.3(c)]{samorodnitsky:sun:2016} for the existence of set $I_A \subsetneq \bbr^d$, for each $A \in \mathscr{R}$.

Let $A \in \mathscr{R}$ be some fixed set. According to \cite{samorodnitsky:sun:2016}, it is said that the distribution $V$ is multivariate subexponential on $A$,  denoted by $V\in\mathcal{S}_A$, if $V_A\in\mathcal{S}$. A key consequence of this definition is the multivariate single big jump principle; see \cite[Cor. 4.10]{samorodnitsky:sun:2016}. In the same manner, by relations \eqref{eq.CCKP.2.10}, \eqref{eq.CCKP.2.14} and \eqref{eq.CCKP.2.15}, if distribution $V_A$ belongs to another heavy-tailed distribution class $\mathcal{B}$, then we say that the distribution $V$ is a corresponding multivariate heavy-tailed distribution on $A$, denoted by $V\in\mathcal{B}_A$, where $\mathcal{B} \in \{\mathcal{C}\cap \mathcal{P_D}\,,\, \mathcal{C}\,,\,\mathcal{D}\cap\mathcal{L}\,,\,\mathcal{D}\,,\,\mathcal{S}\,,\,\mathcal{L}\,,\, \mathcal{P_D}\}$. Please see \cite{konstantinides:passalidis:2024g}, \cite{konstantinides:passalidis:2024h} for the other multivariate heavy-tailed classes. Further we denote $\mathcal{B}_{\mathscr{R}}:=\bigcap_{A \in \mathscr{R}} \mathcal{B}_A$.

Combining \cite[Prop. 2.1]{konstantinides:passalidis:2024g}, \cite[Prop. 3.1]{konstantinides:passalidis:2024h} %
with \cite[Prop. 4.14]{samorodnitsky:sun:2016}, 
together with the univariate inclusions in \eqref{eq.CCKP.2.6}, 
and excluding the class of regularly varying distributions, we obtain
\begin{align*}
{MRV}_{-\alpha}(\mu,G) 
\subsetneq (\mathcal{C}\cap \mathcal{P_D})_{\mathscr{R}}
\subsetneq \mathcal{C}_{\mathscr{R}}
\subsetneq (\mathcal{D}\cap\mathcal{L})_{\mathscr{R}}
\subsetneq \mathcal{S}_{\mathscr{R}} 
\subsetneq \mathcal{L}_{\mathscr{R}},
\end{align*}
for any $\alpha>0$, any distribution $G$, and any limit measure $\mu$.
Moreover, the above inclusions remain valid if $\mathcal{B}_{\mathscr{R}}$ 
is replaced by $\mathcal{B}_A$ for any $A \in \mathscr{R}$ and any 
$\mathcal{B} \in \{\mathcal{C}\cap \mathcal{P_D}\,,\, \mathcal{C}\,,\, 
\mathcal{D}\cap\mathcal{L}\,,\, \mathcal{S}\,,\, \mathcal{L}\}$.
For more properties, examples and applications of these distribution classes, we cite 
\cite{samorodnitsky:sun:2016}, \cite{konstantinides:passalidis:2024g}, \cite{konstantinides:passalidis:2024j}, \cite{konstantinides:liu:passalidis:2025}, \cite{konstantinides:passalidis:2024h}.

\subsection{Dependence modeling} \label{sec.CCKP.2.2}
Here, we provide two dependence structures that can be used to describe the dependency among the claim vectors. As noted earlier, each claim vector is allowed to have dependent components.

Before doing this, we need to introduce the following two dependence structures that describe the dependencies among countably infinitely many random variables. The first dependence structure, introduced by \cite{lehmann:1966}, was found in many applications, see for example \cite{geluk:tang:2009}, \cite{wang:2011}, \cite{jiang:gao:wang:2014},  \cite{geng:liu:wang:2023}.

\begin{definition}
 Let $\{Z_i\,,\;i\geq 1\}$ be a sequence of nonnegative random variables, with distributions $V_{1}\,,\,V_{2}\,,\ldots$, respectively. We say that they are regression dependent, symbolically ${RD}$, if there exist some constants $x_0>0$ and  $K>0$, such that
\beam \label{eq.CCKP.2.17}
\PP\left[ Z_{i} > x_i\;|\;Z_{j}=x_j\,,\;{\text with}\; j \in J\right] \leq K\,\bV_{i}(x_i)\,,
\eeam
holds for any $i\in\mathbb{N}$, $j \in J \subsetneq \mathbb{N}\setminus\{i\}$, with $J \neq \emptyset$, and $x_i\wedge x_j > x_0$.
\end{definition}


The second dependence structure, introduced by \cite{chen:yuen:2009}, was also used in several papers, see for example in \cite{cheng:2014}, \cite{li:2013}, \cite{li:2023b}, \cite{li:2025}, among others.

\begin{definition} 
Let $\{Z_i\,,\;i\geq 1\}$ be a sequence of nonnegative random variables, with distributions $V_{1}\,,\,V_{2}\,,\ldots$, respectively. We say that they are quasi-asymptotically independent, symbolically ${QAI}$, if
\beam \label{eq.CCKP.2.18}
\lim \dfrac{\PP\left[ Z_i > x\,,\;Z_j>x \right]}{\bV_{i}(x)+\bV_{j}(x) }=0\,,
\eeam
holds for any $i,\,j\in\mathbb{N}$, with $i \neq j$. 
\end{definition}

It is easy to check that if random variables $\{Z_i\,,\;i\geq 1\}$ are ${RD}$, then they are also ${QAI}$, but the converse is not generally true. We refer the reader to \cite{cheng:2014} for some counterexamples in this direction. Let us remark that the structures ${RD}$ and ${QAI}$ were defined for random variables with support the whole real line $\bbr$, with small changes with respect to absolute values in relations \eqref{eq.CCKP.2.17} and \eqref{eq.CCKP.2.18}, however, we restrict ourselves to the previous form, for the sake of consistency with the present paper.

Now we are ready to give our definitions of dependence structures among random vectors. 

\begin{definition} \label{ass.CCKP.3.1}
Let $A \in \mathscr{R}$ be a fixed set, $\{{\bf Z}^{(i)}\,,\;i \in \bbn\}$ be a sequence of nonnegative random vectors and $\{{Z}_A^{(i)}\,,\;i \in \bbn\}$ be a sequence of random variables defined by (\ref{eq.CCKP.2.14}). If the $ \{ Z_A^{(i)}\,,\;i \in \bbn \}$ are ${RD}$, then we say that the vectors  $\{{\bf Z}^{(i)}\,,\;i \in \bbn \}$ are regression dependent on $A$, symbolically ${RD}_A$. 
\end{definition} 

Analogously, we can define the following dependence structure based on ${QAI}$.

\begin{definition} \label{ass.CCKP.4.1}
Let $A \in \mathscr{R}$ be a fixed set, $\{{\bf Z}^{(i)}\,,\;i \in \bbn\}$ be a sequence of nonnegative random vectors and $\{{Z}_A^{(i)}\,,\;i \in \bbn\}$ be a sequence of random variables defined by (\ref{eq.CCKP.2.14}). If the $\{ Z_A^{(i)}\,,\;i \in \bbn \} $ are ${QAI}$, then we say that the vectors  $\{{\bf Z}^{(i)}\,,\;i \in \bbn \}$ are quasi-asymptotically independent on $A$, symbolically ${QAI}_A$.
\end{definition}

It is easy to see that if $\{{\bf Z}^{(i)}\,,\;i \in \bbn\}$ are ${RD}_A$, then they are ${QAI}_A$ for the same $A\in\mathscr{R}$. To specifically illustrate the proper inclusion relationship between ${RD}_A$ and ${QAI}_A$, we end this section by giving two tractable examples that satisfy ${QAI}_A$.

\bexam \label{eg1}
For simplicity, we consider the following two-dimensional random vectors ${\bf X}=(X_1\,,\; X_2)$ and ${\bf Y}=(Y_1\,,\; Y_2)$, where each component variable shares the same absolutely continuous marginal distribution $G$. Let 
\beao
A=\left\{(x\,,\;y)\,:\,\, x \vee y >1 \right\}\in\mathscr{R}\,,
\eeao
and the joint distribution of component variables follows the FGM distribution of the form
\begin{align*}
	&\PP\left[X_1>x_1\,,\, X_2>x_2\,,\,Y_1>y_1\,,\,Y_2>y_2\right]\\
	&=\overline{G}(x_1)\,\overline{G}(x_2)\,\overline{G}(y_1)\,\overline{G}(y_2)\left(1+\theta G(x_1)\,G(x_2)\,G(y_1)\,G(y_2)\right)
\end{align*}with $\theta\in[-1,\, 1]$. It is easy to verify that any two and any three of $X_1\,,\;X_2\,,\;Y_1\,,\;Y_2$ are mutually independent. Only for the examples, we define $X_A:=\sup\{u:{\bf X}\in u\,A\}$ and  $Y_A:=\sup\{u:{\bf Y}\in u\,A\}$. Then, for $x,\, y\geq 0$, by Sklar's theorem (see e.g. \cite{nelsen:2006}) one can derive that \begin{align*}
	&\quad \PP\left[X_A>x\,,\; Y_A>y\right]=\PP\left[\left\{X_1>x\right\}\cup \left\{X_2>x\right\}\,,\; \left\{Y_1>y\right\}\cup\left\{Y_2>y\right\}\right]\\[2mm]
	&=4\overline{G}(x)\,\overline{G}(y)-2\overline{G}(x)\,\left(\overline{G}(y)\right)^2-2\left(\overline{G}(x)\right)^2\,\overline{G}(y)+\left(\overline{G}(x)\right)^2\,\left(\overline{G}(y)\right)^2(1+\theta G^2(x)\,G^2(y))\\[2mm]
	&=: \widehat{C}(\overline{G_A}(x),\overline{G_A}(y)),
\end{align*}
where $\widehat{C}(\cdot,\,\cdot)$ denotes the survival copula of $(X_A\,,\, Y_A)$, and 
\begin{align*}
	\overline{G_A}(x)=\PP\left[X_A>x\right]=\PP\left[Y_A>x\right]=2\overline{G}(x)-\left(\overline{G}(x)\right)^2.
\end{align*}

Thus, for  $x$ large enough, we have
\begin{align*}
	&\quad \PP\left[X_A>x| Y_A=y\right]\\
	&=\frac{\partial\widehat{C}(\overline{G_A}(x),\overline{G_A}(y))}{\partial \overline{G_A}(y)}=\frac{\partial\widehat{C}(\overline{G_A}(x),\overline{G_A}(y))}{\partial \overline{G}(y)}\times\frac{\mathrm{d}\overline{G}(y)}{\mathrm{d}\overline{G_A}(y)}\\
	&=\left(2G(y)\right)^{-1}\left(4\overline{G}(x)G(y)-2\left(\overline{G}(x)\right)^2+2\left(\overline{G}(x)\right)^2\overline{G}(y)(1+\theta G^2(x)G^2(y))\right.\\
	&\quad-2\theta \left.\left(\overline{G}(x)\right)^2\left(\overline{G}(y)\right)^2G^2(x)G(y)\right)\\
	&\leq (2\overline{G}(x)+\left(\overline{G}(x)\right)^2)\leq 4\PP\left[X_A>x\right].
\end{align*}
Analogously, for any $x$ large enough, it holds that
\begin{align*}
	\PP\left[Y_A>x| X_A=y\right]\leq 4\PP\left[Y_A>x\right]\,.
\end{align*}
Then, it is obvious that in this example, (\ref{eq.CCKP.2.17}) holds, which implies that ${\bf X}$ and ${\bf Y}$ are ${RD}_A$.
 \eexam

Next, we provide an example where two random vectors are ${QAI}_A$ but not ${RD}_A$.

\bexam \label{eg2.2}
	Consider two two-dimensional positive random vectors ${\bf X}=(X_1\,,\; X_2)$ and ${\bf Y}=(Y_1\,,\; Y_2)$. Assume that the two components of each vector are mutually independent, but for random components sharing the same subscript, the joint distribution is given by
	\begin{align*}
		\PP\left[X_i>x, Y_i>y\right]=\left(1+x^3+y\right)^{-1}\, ,\, i=1\, , \,2\,,
	\end{align*}
	for all $x\geq 0$, $y\geq 0$, and $X_i,\,Y_j$ are independent for each $1 \leq i \neq j \leq 2$. Then, for $i=1,\, 2$, the marginal distributions are 
	\begin{align*}
		\PP\left[X_i>x\right]=\left(1+x^3\right)^{-1},\, P\left[Y_i>y\right]=\left(1+y\right)^{-1},
	\end{align*}
	for all $x\geq 0$, $y\geq 0$\,. Following the argument in Example \ref{eg1}, for all $x\geq 0$, $y\geq 0$, we derive that
	\begin{align*}
		&\PP\left[X_A>x, Y_A>y\right]=2\left(1+x^3+y\right)^{-1}+2\left(1+x^3\right)^{-1}\left(1+y\right)^{-1}\\
		&-2\left(1+x^3+y\right)^{-1}\left(1+y\right)^{-1} -2\left(1+x^3\right)^{-1}\left(1+x^3+y\right)^{-1}+\left(1+x^3+y\right)^{-2}\,, 
	\end{align*}
	\begin{align*}
		\PP\left[X_A>x\right]=2\left(1+x^3\right)^{-1}-\left(1+x^3\right)^{-2}\,,
	\end{align*}
	and 
	\begin{align*}
		\PP\left[Y_A>y\right]=2\left(1+y\right)^{-1}-\left(1+y\right)^{-2}\,.
	\end{align*}
	Thus, for $x$ large enough, setting y=x, it holds that 
	\begin{align*}
		\frac{\PP\left[X_A>x, Y_A>x\right]}{\PP\left[X_A>x\right]+\PP\left[Y_A>x\right]}=O(1)\frac{\left(1+x^3+x\right)^{-1}+\left(1+x^3\right)^{-1}\left(1+x\right)^{-1}}{\left(1+x^3\right)^{-1}+\left(1+x\right)^{-1}}\to 0\,,
	\end{align*}
	which implies that ${\bf X}=(X_1\,,\; X_2)$ and ${\bf Y}=(Y_1\,,\; Y_2)$ are ${QAI}_A$. However, they are not ${RD}_A$. Indeed, for $x$ large enough, we have
	\begin{align*}
		\frac{\PP\left[X_A>x, Y_A>x\right]}{\PP\left[X_A>x\right]}\geq \frac{3}{4}\cdot\frac{\left(1+x^3+x\right)^{-1}+\left(1+x^3\right)^{-1}\left(1+x\right)^{-1}}{\left(1+x^3\right)^{-1}}\to \frac{3}{4}>0.
	\end{align*}
	This contradicts \eqref{eq.CCKP.2.17}.
\eexam

\section{Main results}

In this section, we present the main results that establish relation \eqref{eq.CCKP.1.3} under different conditions. 
In what follows, for some fixed $A \in \mathscr{R}$, the claim vectors ${\bf X}^{(i)}$ are related with the random variables
\begin{align*}
 Y_A^{(i)}:= \sup\left\{ u\;:\;{\bf X}^{(i)} \in u\,A \right\}\,,
\end{align*}
that follow the distributions $F_A^{(i)} \equiv F_A$, respectively.

\subsection{On finite time horizon asymptotics with multivariate subexponential claims}
This subsection will establish \eqref{eq.CCKP.1.3} when $F$ is multivariate subexponential and $\{\xi(t)\,,\;t\geq 0\}$ has non-decreasing sample paths. Note that we do not have any additional conditions on the components of each claim vector, which, in combination with ${RD}_A$, allows interdependence among claim vectors.

\bth \label{th.CCKP.3.1}
Consider the aggregate claims from relation \eqref{eq.CCKP.1.2}. Let $ A\in\mathscr {R}$ be a fixed set. We suppose that Assumption \ref{ass.CCKP.1.1} is true, $\{{\bf X}^{(i)}\,,\;i \in \bbn\}$ are ${RD}_A$, $F \in \mathcal{S}_A$, and $\{\xi(t)\,,\;t\geq 0\}$ has non-decreasing sample paths. Then for any fixed $T \in \varLambda\setminus \{\infty\}$, relation \eqref{eq.CCKP.1.3} holds.
\ethe 

\bre \label{rem.CCKP.3.1}
We observe that some similar results, as in relation \eqref{eq.CCKP.1.3} are found in \cite[Th. 5.1]{konstantinides:passalidis:2024g}, where the  relation \eqref{eq.CCKP.1.3} was established in a model with independent claim vectors and a common Poisson counting process, but there was no restriction for non-decreasing sample paths for $\{\xi(t)\,,\;t\geq 0\}$. Furthermore, in \cite[Th. 3.1]{konstantinides:passalidis:2024j} relation \eqref{eq.CCKP.1.3} appears under a larger dependence structure between the claim vectors than ${RD}_A$, but the distribution class of claim vectors is restricted to $ (\mathcal{D}\cap\mathcal{L})_A$. 
\ere

\bre \label{rem.CCKP.3.2}
In Theorem \ref{th.CCKP.3.1}, if we replace the assumption $F \in \mathcal{S}_A$ by $F \in \mathcal{S}_{\mathscr{R}}$, then the conclusion ${RD}_A$ holds for every $A \in \mathscr{R}$, and relation \eqref{eq.CCKP.1.3} remains valid for all $A \in \mathscr{R}$. Similar statements hold for all subsequent main results. For brevity, we shall not reiterate throughout the paper.
\ere
\subsection{On finite time horizon asymptotics with multivariate consistently varying claims}
The following theorem revisits relation \eqref{eq.CCKP.1.3} under a different dependence structure, namely when the claim vectors are ${QAI}_A$ and $F$ belongs to the distribution class $\mathcal{C}_A$, for some fixed $A\in\mathscr{R}$. Compared with the results of the previous section, we see that the claim vector distribution class is reduced from $\mathcal{S}_A$ into $\mathcal{C}_A$, while the dependence structure is generalized from ${RD}_A$ to ${QAI}_A$. Moreover, the condition of non-decreasing sample paths for $\{\xi(t),\, t\ge 0\}$ is removed.
\bth \label{th.CCKP.4.1}
Consider the aggregate claims from relation \eqref{eq.CCKP.1.2}. Let $ A\in\mathscr {R}$ be a fixed set. We suppose that Assumptions \ref{ass.CCKP.1.1} and \ref{ass.CCKP.1.2} are true,  $\{{\bf X}^{(i)}\,,\; i \in \bbn\}$ are ${QAI}_A$, and $F \in \mathcal{C}_A$.
Then for any fixed $T \in \varLambda\setminus\{\infty\}$, relation \eqref{eq.CCKP.1.3} holds.
\ethe 

The following corollary is established by Theorem \ref{th.CCKP.4.1}, building on the fact that $F \in {MRV}_{-\alpha}(\mu,\,G)$. As mentioned previously, it ensures $F_A \in \mathcal{R}_{-\alpha}$, for any $A \in \mathscr{R}$. The dominated convergence theorem and Breiman's theorem, see in \cite{breiman:1965}, become applicable through Assumption \ref{ass.CCKP.1.2}. Restricting ourselves in ${MRV}$, we gain a more explicit expression in relation \eqref{eq.CCKP.1.3}.

\bco \label{cor.CCKP.4.1}
Under the conditions of Theorem  \ref{th.CCKP.4.1}, if $F \in {MRV}_{-\alpha}(\mu,\,G)$, then for any fixed $T \in \varLambda\setminus \{\infty\}$, we obtain the asymptotic relation
\beam \label{eq.CCKP.4c.1}
\PP\left[ {\bf D}(T) \in x\,A  \right] \sim \mu(A)\overline{G}(x)\int_0^T \E\left[ e^{-\alpha\,\xi(s)}\right]\,\lambda(ds)\,. 
\eeam
\eco

\begin{remark}
We note that relation \eqref{eq.CCKP.4c.1} is indeed a precise asymptotic estimate. By \cite[Prop. 4.14]{samorodnitsky:sun:2016}, if $F\in {MRV}_{-\alpha}(\mu,G)$, then, for any $A\in\mathscr R$,
\[
\overline F_A(x)=\mathbb P(Y_A>x)=\mathbb P({\bf X}\in xA)\sim
\mu(A)\overline G(x)\,,
\]
with $0<\mu(A)<\infty$. Hence $F_A\in\mathcal R_{-\alpha}$. This observation justifies the use of Breiman’s theorem in the proof of Corollary \ref{cor.CCKP.4.1}. In particular, the conclusion remains valid for the rare sets $A\in\mathscr R$, even when the components of $\,{\bf X}$ are asymptotically independent.
\end{remark}
\subsection{On infinite time horizon asymptotics}
In this subsection, we deal with the estimate \eqref{eq.CCKP.1.3} over the infinite time horizon, namely for $T=\infty$. Let $A\in\mathscr{R}$ be a fixed set. Assuming that $\{{\bf X}^{(i)}\,,\;i\in \bbn\}$ still represents ${QAI}_A$ sequence of vectors, but the class of the common distribution $F$, is restricted from $\mathcal{C}_A$ into $(\mathcal{C} \cap \mathcal{P_D})_A$. However, we replace Assumption \ref{ass.CCKP.1.2} by Assumption \ref{ass.CCKP.5.1} presented below, which seems rather general to provide more flexibility with respect to the insurer's investment portfolio.

We begin with the basic assumption for the general process $\{\xi(t)\,,\;t \geq 0\}$. Let us remind that $F_A$ has $J_{F_A}^- >0$, since $F  \in \mathcal{P_D}_A$.

\begin{assumption} \label{ass.CCKP.5.1}
Let $\{\xi(t)\,,\;t \geq 0\}$ be a general process. We assume that there exist some $p_1,\,p_2 $ with $0<p_1 < J_{F_A}^- \leq J_{F_A}^+ < p_2 <\infty $, such that
\begin{align}\label{eq.CCKP.5.2} 
	\sum_{i=1}^{\infty} \left(\E\left[ e^{-p_1\,\xi(\tau_i)} \right]\vee \E\left(e^{-p_2\,\xi(\tau_i)} \right] \right)^{1/\rho} < \infty\,,
\end{align}
where $\rho={\bf 1}_{\left(0<J_{F_A}^+<1\right)}+ p_2{\bf 1}_{\left(J_{F_A}^+\geq 1\right)}.$
\end{assumption} 
Thus, through Assumption \ref{ass.CCKP.5.1}, we get rid of the condition of lower bound for $\{\xi(t)\,,\;t \geq 0\}$. It is easy to see that if $\xi(t)=r\,t$, with $r> 0$ the constant interest force, then Assumption \ref{ass.CCKP.5.1} is satisfied. In general, the removal of Assumption \ref{ass.CCKP.1.2} permits the insurer to make riskier investments, which is sometimes necessary to make the insurance company reliable. Especially in cases where insurers have portfolios with heavy-tailed claims, they are encouraged to choose riskier investments.

Now, we provide a general example of stochastic processes that satisfy Assumption \ref{ass.CCKP.5.1}. Indeed, in this example we find out that all L\'{e}vy processes, with Laplace exponent $\phi(p_2)<0$ for some $p_2 > J_{F_A}^+$, satisfy Assumption \ref{ass.CCKP.5.1}. 

\bexam \label{exam.CCKP.5.1}
Let $\{\xi(t)\,,\; t\geq 0\}$ be a real-valued L\'{e}vy process. The Laplace exponent is defined by the relation
\beao
\phi(k)=\log \E\left[ e^{-k\,\xi(1)} \right]\,,
\eeao
for any $k \in \bbr$. It is known that if $\phi(k) < \infty$, then 
\beam  \label{eq.CCKP.5.1e}
\E\left[ e^{-k\,\xi(t)} \right]= e^{t\,\phi(k)} < \infty\,,
\eeam
see for example in \cite[Prop. 3.14]{cont:tankov:2004}.

Let $\phi(p_2) < 0$, for some $p_2 > J_{F_A}^+$. Then, from relation \eqref{eq.CCKP.5.1e}, and taking into account that the $\{\tau_i\,,\;i \in \bbn\}$ represent a renewal counting process, we obtain that it holds
\begin{align*} 
\E\left[ e^{-p_2\,\xi(\tau_i)} \right]=\left( \E\left[ e^{\tau_1\,\phi(p_2)} \right] \right)^i < 1,
\end{align*}
for any $i \in \bbn$. Further, for $p_1 \in (0,\,J_{F_A}^-)$, it holds
\begin{align*}  
\E\left[ e^{-p_1\,\xi(\tau_i)} \right]=\left( \E\left[ e^{\tau_1\,\phi(p_1)} \right] \right)^i < 1,
\end{align*}
where in the last step we used the fact that $\phi(p_1)< 0$, since $\phi(p_2)< 0$, $\phi(0)= 0$ and the $\phi(k)$ is convex with respect to $k$, that implies $\phi(k)<0$, for any $k \in (0,\,p_2]$.

Then, for $0<p_1 < J_{F_A}^- \leq J_{F_A}^+ < p_2 <\infty $ and $\rho={\bf 1}_{\left(0<J_{F_A}^+<1\right)}+ p_2{\bf 1}_{\left(J_{F_A}^+\geq 1\right)}$, it is not hard to verify that 
\begin{align*}
	&\sum_{i=1}^{\infty} \left(\E\left[ e^{-p_1\,\xi(\tau_i)} \right]\vee \E\left(e^{-p_2\,\xi(\tau_i)} \right] \right)^{1/\rho} \leq \sum_{i=1}^{\infty} \left[\left(\E\left[  e^{\tau_1\,\phi(p_1)} \right] \right)^i + \left(\E\left[  e^{\tau_1\,\phi(p_2)}  \right] \right)^i\right]^{1/\rho}\\[2mm] \notag
&\leq \sum_{i=1}^{\infty} \left[\left(\E\left[  e^{\tau_1\,\phi(p_1)} \right] \right)^{i/\rho} + \left(\E\left[  e^{\tau_1\,\phi(p_2)}  \right] \right)^{i/\rho}\right]< \infty\,,
\end{align*}
where in the penultimate step the elementary inequality $(a+b)^r \leq a^r + b^r$, for any $a,\,b \geq 0$, and any $r \in (0,\,1]$. Which indicates that Assumption \ref{ass.CCKP.5.1} holds.
\eexam

\bre \label{rem.CCKP.5.0}
We notice that the assumption that the logarithmic prices of the investment portfolio is L\'{e}vy processes, with Laplace exponent that satisfies the condition $\phi(p)<0$, for some $p$ greater than the upper Matuszewska index of distribution of the claims, is popular enough in one-dimensional and multidimensional risk models, that study cases, where the insurers make risky investments, see for example in \cite{tang:wang:yuen:2010}, \cite{yang:wang:konstantinides:2014}, \cite{li:2016}, \cite{yang:su:2023} among others. Therefore, through Example \ref{exam.CCKP.5.1}, we realize the generality of Assumption \ref{ass.CCKP.5.1}.  
\ere

Now, we are ready for the main result of this subsection.

\bth \label{th.CCKP.5.1}
Consider the aggregate claim vector \eqref{eq.CCKP.1.2} with $T=\infty$. Let $A \in \mathscr{R}$ be a fixed set. We supposed that Assumptions \ref{ass.CCKP.1.1} and \ref{ass.CCKP.5.1} are true, $\{{\bf X}^{(i)}\,,\;i\in \bbn\}$ are ${QAI}_A$, and $F \in (\mathcal{C} \cap \mathcal{P_D})_A$. Then \eqref{eq.CCKP.1.3} holds with $T=\infty$. 
\ethe

If Assumption \ref{ass.CCKP.5.1} holds, then it is obvious that
\beao
\E\left[ e^{-p_1\,\xi(\tau_i)} \right]\vee \E\left[e^{-p_2\,\xi(\tau_i)} \right] \leq \left(\E\left[ e^{-p_1\,\xi(\tau_i)} \right]\vee \E\left[e^{-p_2\,\xi(\tau_i)} \right] \right)^{1/\rho}< \infty\,,
\eeao 
hence if $F_A \in \mathcal{R}_{-\alpha}$, for some $\alpha \in (0,\,\infty)$, then it holds that
\beao
 \E\left[e^{-\alpha\,\xi(\tau_i)} \right]\leq \E\left[ e^{-p_1\,\xi(\tau_i)} \right]\vee \E\left[e^{-p_2\,\xi(\tau_i)} \right] < \infty\,,
\eeao 
for some $0<p_1<\alpha<p_2<\infty$. Therefore, we can analogously apply the dominated convergence theorem and Breiman's theorem in Theorem \ref{th.CCKP.5.1} to establish the following corollary.

\bco \label{cor.CCKP.5.1}
Under the conditions of  Theorem \ref{th.CCKP.5.1}, if $F \in {MRV}_{-\alpha}(\mu,\,G)$, then it holds that
\begin{align*}
{\PP[{\bf D}(\infty) \in x\,A]}\sim\mu({A}){\overline{G}(x)}\int_{0}^{\infty} \E\left[e^{-\alpha\,\xi(s)} \right]\,\lambda(ds)\,.
\end{align*}
\eco

\bre \label{rem.CCKP.5.a^*}
Relation \eqref{eq.CCKP.4c.1} for $T<\infty$ and $T=\infty$, it can take more explicit forms in the case where 
the $A \in\mathscr{R}$, the dependence structure, and the distribution class of the components of 
vector ${\bf X}$ are determined. See in \cite[Sec. 4]{konstantinides:liu:passalidis:2025} for some examples in this direction. 
 The class ${MRV}$ gives extremely explicit relations (as we see in Corollaries \ref{cor.CCKP.4.1} and \ref{cor.CCKP.5.1}), 
since the dependence structures are depicted only through the Radon measure $\mu$, and the speed of convergence of the tail is given only by 
a distribution $G \in \mathcal{R}_{-\alpha}$. However, the wider frame of our main results has several advantages. For 
example the ${MRV}$ requires in some sense ``equivalent tails" of marginals, while the class $\mathcal{C}_A$ (as also the class $(\mathcal{C} \cap \mathcal{P_D})_A$)
 does not.
Let us consider, as an example, the set $A_2$ of \eqref{eq.CCKP.2.12} with $d=2$, $l_1=l_2=1/2$, and $X_1,\,X_2$ non-negative 
arbitrarily dependent random variables with distributions $F_1,\,F_2$, such that $F_1 \in \mathcal{C}$ (or $F_1 
\in \mathcal{C} \cap \mathcal{P_D}$) and $\bF_2(x) =o[\bF_1(x)]$, then following similar arguments as that in \cite[Exam. 4.4]{konstantinides:liu:passalidis:2025}, we obtain that 
\beao
\bF_{A_2}(x) = \PP\left( \dfrac 12\,X_1 + \dfrac 12\,X_2 >x \right) \sim \PP\left( \dfrac 12\,X_1 >x \right)\,,
\eeao   
hence $F_{A_2} \in \mathcal{C}$ (or $F_{A_2} \in \mathcal{C} \cap \mathcal{P_D}$, respectively) and consequently we find 
$F \in \mathcal{C}_{A_2}$ (or $F \in (\mathcal{C} \cap \mathcal{P_D})_{A_2}$, respectively). So, the classes $\mathcal{C}_{A_2}$ or $(\mathcal{C} \cap \mathcal{P_D})_{A_2}$, 
that were used in Theorem \ref{th.CCKP.4.1} and Theorem \ref{th.CCKP.5.1} do not exclude  subexponential marginals and 
cases of non-equivalent tails of the marginals, as it happens with ${MRV}$. 
We can see that in the case of $F_{A_2}$, if $F_1 \in  \mathcal{R}_{-\alpha}$, for $\alpha \in (0,\,\infty)$, let say Pareto distribution, 
and $F_2 \in \mathcal{S} \setminus \mathcal{R}$, e.g., the lognormal or the Weibull distribution, then we have 
$F_{A_2} \in  \mathcal{R}_{-\alpha} \subsetneq \mathcal{C}$ but $F \notin {MRV}_{-\alpha}(\mu,\,G)$. We note also that in this last example $F$ does not belong even in the nonstandard ${MRV}$ class, see \cite{chen:yang:2019} and \cite{tang:yang:2019} for more details on this class.
\ere

\bre \label{rem.CCKP.5.3}
The first entrance time of one-dimensional and multidimensional random walks into rare sets has been well studied and is directly related to the insurer's surplus process; see \cite{liu:woo:2014} and the references therein.
However, due to the fact that the study is restricted to random walks only, it permits covering several situations in the modern insurance industry, such as the presence of investments and interdependence. Therefore, our results present a first picture for the conditional distribution of the first entrance time of the aggregate claims into some rare set $x\,A \in \mathscr{R}$.

Firstly, for some $A \in \mathscr{R}$, we define by
\beam \label{eq.CCKP.5.c1} 
\tau(x) := \inf \left\{ t>0\;:\; {\bf D}(t) \in x\,A \right\}\,,
\eeam
the first entrance time, which could eventually be a defective random variable, since we could have that 
\begin{align*}
\PP[\tau(x) < \infty]=\PP[{\bf D}(\infty) \in x\,A]<1\,,
\end{align*}
where at the first step we used relation \eqref{eq.CCKP.5.c1}, the increasing property of set $A$ and the nonnegativeness of the summands of ${\bf D}(\infty)$. Thus, in case of simple interest, when $\xi(t)=r\,t$, with $r>0$ and under the condition $F \in (\mathcal{C} \cap \mathcal{P_D})_A$, that implies the validity of Theorems \ref{th.CCKP.4.1} and \ref{th.CCKP.5.1}, we obtain that
\beao
&& \PP\left[ \tau(x) \leq T\;|\; \tau(x) < \infty \right] = \dfrac{\PP\left[ \tau(x) \leq T \right] }{\PP\left[  \tau(x) < \infty \right] }= \dfrac{\PP\left[ {\bf D}(T) \in x\,A \right] }{\PP\left[ {\bf D}(\infty) \in x\,A\right] } \\[2mm]
&&\sim \dfrac{\int_{0}^{T} \PP\left[ {\bf X}\,e^{-r\,s}  \in x\,A \right]\,\lambda(ds) }{\int_{0}^{\infty} \PP\left[ {\bf X}\,e^{-r\,s}  \in x\,A \right]\,\lambda(ds)} 
\eeao
for any $T\in \Lambda \setminus\{\infty\}$.

Furthermore, if we assume that the counting process $\{N(t)\,,\;t\geq 0\}$ is Poisson, and  under the condition $F_A\in\mathcal{R}_{-\alpha}$, for some  $\alpha \in (0,\,\infty)$ (more general than condition $F\in {MRV}_{-\alpha}(\mu,\,G)$), then we obtain 
\beam \label{eq.CCKP.5.c3}
\PP[\tau(x) \leq T\;|\;\tau(x) < \infty ] \sim \dfrac{\int_0^T e^{-\alpha\,r\,s} ds}{\int_0^{\infty} e^{-\alpha\,r\,s} ds}=1 -e^{-\alpha\,r\,T}\,%
\eeam 
for any $T > 0$. In the first step, we used the dominated convergence theorem via the uniformity of $\mathcal{R}_{-\alpha}$. That means, in case of constant interest force, with $F_A\in\mathcal{R}_{-\alpha}$, and Poisson counting process $\{N(t)\,,\;t\geq 0\}$, we find by relation \eqref{eq.CCKP.5.c3} that the limit distribution of $\tau(x)$ conditioned by $\tau(x)<\infty$ coincides with the exponential one with parameter $\alpha\,r$. 
\ere

\bre \label{rem.CCKP.5.4}
The sets $A \in \mathscr{R}$ are directly related to some ``ruin"-sets, see \cite[Ass. 5.1]{samorodnitsky:sun:2016}. Hence, in the case where the premium densities are bounded from above, e.g., linear premiums, the results of Theorems \ref{th.CCKP.3.1} and \ref{th.CCKP.4.1} can be straightforwardly translated into the finite time ruin probabilities, through some classical methods; see, for example, in the proof of \cite[Cor. 3.1]{konstantinides:liu:passalidis:2025} or of \cite[Cor. 3.1]{konstantinides:passalidis:2024j}. For the case of Theorem \ref{th.CCKP.5.1}, the corresponding infinite time ruin probabilities do not follow directly, in the general case of Assumption \ref{ass.CCKP.5.1}. However, in special cases, where the $\{\xi(t)\,,\;t\geq 0 \}$ satisfies simultaneously Assumption \ref{ass.CCKP.1.2}, for any $T>0$ and Assumption \ref{ass.CCKP.5.1} (e.g., the $\{\xi(t)\,,\;t\geq 0 \}$ represents a non-negative L\'{e}vy process), then the infinite time ruin probability can be established relatively easily. 
\ere

\section{Proofs of the main results}
We now proceed to the proofs of our main results. However, before proving each theorem, we require some auxiliary lemmas.

\subsection{Proof of Theorem \ref{th.CCKP.3.1}} \label{sec.CCKP.3.2}

We begin with some preliminary lemmas. The first one has its own merit, as it demonstrates the insensitivity property of the multivariate linear single big jump principle of the scale mixture sums with respect to ${RD}_A$. We note also that the Lemma \ref{lem.CCKP.3.1} adopts more general assumptions than what we need in the proof of Theorem \ref{th.CCKP.3.1}.

\ble \label{lem.CCKP.3.1}
Let $A \in \mathscr{R}$ be a fixed set. We suppose that the random vectors ${\bf Z}^{(1)},\,\ldots,\,{\bf Z}^{(n)}$ are ${RD}_A$, with distributions $V_1,\,\ldots,\,V_n \in \mathcal{L}_A$, respectively. Let $W_1,\,\ldots,\,W_n$ be $n$ arbitrarily dependent random variables (nondegenerate to zero), which are independent of ${\bf Z}^{(1)},\,\ldots,\,{\bf Z}^{(n)}$, such that $0\leq W_i\leq b$ almost surely for some $0< b < \infty$ and for any $i=1,\,\ldots,\,n$. If additionally we assume that $V_i(xA)\asymp V(xA)$, for some $V\in\mathcal{S}_A$, and for any $i=1,\,\ldots,\,n$ then it holds that
\beam \label{eq.CCKP.3.3}
\PP\left[ \sum_{i=1}^n W_i\,{\bf Z}^{(i)} \in x\,A \right] \sim \sum_{i=1}^n\PP\left[  W_i\,{\bf Z}^{(i)} \in x\,A \right] \,.
\eeam
\ele

\noindent\pr~
Relation (\ref{eq.CCKP.3.3}) amounts to the conjunction of
\beam \label{eq.CZT.3.1}
\PP\left[ \sum_{i=1}^n W_i\,{\bf Z}^{(i)} \in x\,A \right] \lesssim \sum_{i=1}^n\PP\left[  W_i\,{\bf Z}^{(i)} \in x\,A \right] \,,
\eeam
and 
\beam \label{eq.CZT.3.2}
\PP\left[ \sum_{i=1}^n W_i\,{\bf Z}^{(i)} \in x\,A \right] \gtrsim \sum_{i=1}^n\PP\left[  W_i\,{\bf Z}^{(i)} \in x\,A \right] \,.
\eeam
We prove (\ref{eq.CZT.3.1}) in two steps, first we assume that $W_1,\,\ldots,\,W_n$ are positive, i.e., 
\beao
\PP(0<W_i\leq b)=1\,,
\eeao 
for any $i=1,\,\ldots,\,n$. For any subsets $I\subseteq \mathbb{I}=\{1,2,\ldots,n\}$ and any $0<\varepsilon<1$, write 
\begin{align*}
	\Omega_{I}^{\varepsilon}(W)=\{\omega:W_i(\omega)>\varepsilon\,\, {\rm for}\,\, i\in I,\, {\rm and}\,\, W_j(\omega)\leq \varepsilon\,\,{\rm for}\,\, j\notin I\}.
\end{align*}
Then, for the (\ref{eq.CZT.3.1}), we have by following the lines of the proof of \cite[Th. 1 ]{tang:yuan:2014} that
\begin{align*}
	&\quad \PP\left[ \sum_{i=1}^n W_i\,{\bf Z}^{(i)} \in x\,A \right] \leq \PP\left[\sum_{i=1}^n W_i\,Z_A^{(i)} > x \right]\\[2mm]
	&\leq\sum_{I\subseteq\mathbb{I}}\PP\left[\sum_{i=1}^n W_i\,Z_A^{(i)} > x ,\,\Omega_{I}^{\varepsilon}(W)\right]\leq \sum_{I\subseteq\mathbb{I}}\PP\left[\sum_{i\in I} W_i\,Z_A^{(i)} +\sum_{j\notin I} \varepsilon\, Z_A^{(j)} > x,\, \Omega_{I}^{\varepsilon}(W)\right]\nonumber\\[2mm]
	&\sim \sum_{I\subseteq\mathbb{I}}\sum_{i\in I} \PP\left[ W_i\,Z_A^{(i)}>x,\, \Omega_{I}^{\varepsilon}(W)\right] +\sum_{I\subseteq\mathbb{I}}\sum_{j\notin I} \PP\left[\varepsilon\,Z_A^{(j)} >x,\, \Omega_{I}^{\varepsilon}(W)\right]\nonumber\\[2mm]
	&=\sum_{i=1}^n\sum_{I\subseteq \mathbb{I}:\,i\in I}\PP\left[ W_i\,Z_A^{(i)}>x,\, \Omega_{I}^{\varepsilon}(W)\right]+\sum_{j=1}^n\sum_{I\subseteq \mathbb{I}:\,j\notin I}\PP\left[\varepsilon\,Z_A^{(j)} >x,\, \Omega_{I}^{\varepsilon}(W)\right]\nonumber\\[2mm]
	&  =\sum_{i=1}^n\sum_{I\subseteq \mathbb{I}:\,i\in I}\PP\left[ W_i\,Z_A^{(i)}>x,\, \Omega_{I}^{\varepsilon}(W)\right]+\sum_{j=1}^n\sum_{I\subseteq \mathbb{I}:\,j\notin I}\PP\left[\varepsilon\,Z_A^{(j)} >x,\, W_j>\vep \right]\dfrac{\PP[\Omega_{I}^{\varepsilon}(W)]}{\PP[W_j>\vep]} \nonumber\\[2mm]
	&\leq \left(1+\max_{1\leq i\leq n }\frac{\PP[W_i\leq \varepsilon]}{\PP[W_i>\varepsilon]}\right)\sum_{i=1}^n\PP\left[ W_i\,Z_A^{(i)} > x \right]\nonumber\\[2mm]
	&=\left(1+\max_{1\leq i\leq n }\frac{\PP[W_i\leq \varepsilon]}{\PP[W_i>\varepsilon]}\right)\sum_{i=1}^n\PP\left[ W_i\,{\bf Z}^{(i)} \in x\,A \right],
\end{align*}
where in the first step, we used \cite[Prop. 2.4]{konstantinides:passalidis:2024g}, and in the fourth step, we used \cite[Th. 3.2]{geng:liu:wang:2022}.
By letting $\varepsilon\downarrow 0$, we conclude (\ref{eq.CZT.3.1}). 

In the second step, we prove the case where \(W_1,\ldots,W_n\) may take
the value zero with positive probability. For any subset
\(I\subseteq\mathbb I:=\{1,\ldots,n\}\), write
\[
\Omega_I^0(W)
=
\{\omega: W_i(\omega)>0 \ {\rm for}\ i\in I,\ 
W_j(\omega)=0 \ {\rm for}\ j\in \mathbb I\setminus I\}.
\]
The events \(\Omega_I^0(W)\), \(I\subseteq\mathbb I\), form a partition of
the sample space. On \(\Omega_\emptyset^0(W)\), all weights are equal to zero,
and this event gives no contribution. Hence,
\[
\begin{aligned}
\mathbb P\left[
\sum_{i=1}^n W_i{\bf Z}^{(i)}\in xA
\right]
\le
\mathbb P\left[
\sum_{i=1}^n W_iZ_A^{(i)}>x
\right] =
\sum_{\emptyset\neq I\subseteq\mathbb I}
\mathbb P\left[
\sum_{i\in I}W_iZ_A^{(i)}>x,\ \Omega_I^0(W)
\right].
\end{aligned}
\]
Fix a nonempty \(I\subseteq\mathbb I\). If
\(\mathbb P(\Omega_I^0(W))=0\), then the corresponding term is zero. Otherwise,
define
\[
\mathbb P_I(\cdot):=\mathbb P(\cdot\mid \Omega_I^0(W)).
\]
Since \(\Omega_I^0(W)\) is determined only by the weights, that are
independent of \({\bf Z}^{(1)},\ldots,{\bf Z}^{(n)}\), under \(\mathbb P_I\)
the distribution of the random vectors \({\bf Z}^{(1)},\ldots,{\bf Z}^{(n)}\)
is unchanged. In particular, their \(RD_A\) structure is preserved, and the
active weights \(\{W_i:i\in I\}\) remain independent of the random vectors.
Moreover, under \(\mathbb P_I\),
\[
0<W_i\le b,\qquad i\in I,
\]
almost surely. Therefore, the positive-weight case proved above can be applied
to the active family \(\{W_i{\bf Z}^{(i)}:i\in I\}\) under \(\mathbb P_I\), and
we obtain
\[
\mathbb P_I\left[
\sum_{i\in I}W_iZ_A^{(i)}>x
\right]
\sim
\sum_{i\in I}
\mathbb P_I\left[
W_iZ_A^{(i)}>x
\right].
\]
Multiplying both sides by \(\mathbb P(\Omega_I^0(W))\), it follows that
\[
\mathbb P\left[
\sum_{i\in I}W_iZ_A^{(i)}>x,\ \Omega_I^0(W)
\right]
\sim
\sum_{i\in I}
\mathbb P\left[
W_iZ_A^{(i)}>x,\ \Omega_I^0(W)
\right].
\]
Since there are only finitely many nonempty subsets \(I\subseteq\mathbb I\),
summing over all such \(I\) yields
\[
\begin{aligned}
&\mathbb P\left[
\sum_{i=1}^n W_i{\bf Z}^{(i)}\in xA
\right]
\lesssim
\sum_{\emptyset\neq I\subseteq\mathbb I}
\sum_{i\in I}
\mathbb P\left[
W_iZ_A^{(i)}>x,\ \Omega_I^0(W)
\right] \\
&=
\sum_{i=1}^n
\sum_{\substack{\emptyset\neq I\subseteq\mathbb I\\ i\in I}}
\mathbb P\left[
W_iZ_A^{(i)}>x,\ \Omega_I^0(W)
\right]=
\sum_{i=1}^n
\mathbb P\left[
W_iZ_A^{(i)}>x,\ W_i>0
\right] \\
&=
\sum_{i=1}^n
\mathbb P\left[
W_iZ_A^{(i)}>x
\right] =
\sum_{i=1}^n
\mathbb P\left[
W_i{\bf Z}^{(i)}\in xA
\right].
\end{aligned}
\]
This proves the upper bound in the case where some weights may be zero.
This completes the proof of (\ref{eq.CZT.3.1}). As for (\ref{eq.CZT.3.2}), because $x\,A$ is an increasing set, $W_i$, ${\bf Z}^{(i)} $ are nonnegative, and by Bonferroni's inequality, we find that it holds
\begin{align} \label{eq.CCKP.3.5}
&\PP\left[ \sum_{i=1}^n W_i\,{\bf Z}^{(i)} \in x\,A \right] \geq\, \PP\left[ \bigcup_{i=1}^n \{ W_i\,{\bf Z}^{(i)} \in x\,A \} \right] \geq\\[2mm]  \nonumber
&\sum_{i=1}^n\PP\left[  W_i\,{\bf Z}^{(i)} \in x\,A \right]- \sum_{1\leq i<j\leq n} \PP\left[  W_i\,{\bf Z}^{(i)} \in x\,A\,,\;W_j\,{\bf Z}^{(j)} \in x\,A \right] \sim \sum_{i=1}^n\PP\left[  W_i\,{\bf Z}^{(i)} \in x\,A \right]\,,
\end{align}
where the last step is justified, since each term of the second sum (in the penultimate step) is asymptotically negligible with respect to each term of the first sum. Indeed, for any $1 \leq i \neq j \leq n$ and for any $x/b > x_0$, by the assumption that ${\bf Z}^{(i)},\,{\bf Z}^{(j)} $ are ${RD}_A$, it follows that
\beao
&&\PP\left[  W_i\,{\bf Z}^{(i)} \in x\,A\,,\;W_j\,{\bf Z}^{(j)} \in x\,A \right]=\PP\left[  W_i\,Z_A^{(i)} > x\,,\;W_j\,Z_A^{(j)} > x \right]\\[2mm] 
&& \leq \PP[b\,Z_A^{(j)} > x]\,\PP\left[  W_i\,Z_A^{(i)} > x\;|\;b\,Z_A^{(j)} > x \right] \\[2mm] 
&& =\PP[b\,Z_A^{(j)} > x]\,\int_0^b \PP\left[ y\,Z_A^{(i)} > x\;|\;b\,Z_A^{(j)} > x \right] \PP[W_i \in dy]\\[2mm] 
&& =O(1)\, \PP[b\,Z_A^{(j)} > x]\,\int_0^b \PP\left[ y\,Z_A^{(i)} > x \right] \,\PP[W_i \in dy] \\[2mm]
&& =O(1)\, \PP[b\,Z_A^{(j)} > x]\,\PP\left[ W_i\,Z_A^{(i)} > x \right]  =o\left( \PP\left[ W_i\,Z_A^{(i)} > x \right] \right) = o\left( \PP\left[  W_i\,{\bf Z}^{(i)} \in x\,A \right] \right)\,.
\eeao
Hence, we obtain relation (\ref{eq.CZT.3.2}) and close the proof of this lemma. 
~\halmos


The following lemma is a restatement of the generalized Kesten's inequality from \cite[Th. 1]{geng:liu:wang:2023}. We provide this interpretation under stricter conditions to fit better with the needs of our paper.

\ble \label{lem.CCKP.3.2} 
Let $\{Z_{i}\,,\;i\in \bbn\}$ be ${RD}$ and nonnegative random variables with common distribution $V \in \mathcal{S}$. Let $\{W_i\,,\;i\in \bbn\}$ be non-negative, nondegenerate to zero, and uniformly upper-bounded random weights that are independent of $\{Z_{i}\,,\;i\in \bbn\}$. Then for any $\vep > 0$, there exists some constant $C=C(\vep)> 0$, such that 
\begin{align*}
\PP\left[ \sum_{i=1}^n W_i\, Z_{i} >x \right] \leq C\,(1+\vep)^n\,\sum_{i=1}^n \PP\left[ W_i\, Z_{i} > x \right] \,,
\end{align*}
holds for all $n\in \bbn$ and all $x \geq 0$.
\ele 

Now we proceed to the proof of Theorem \ref{th.CCKP.3.1}.

\noindent{\bf Proof of Theorem \ref{th.CCKP.3.1}}~
Let $M$ be some large integer. From the total probability theorem, we obtain the equality
\begin{align} \label{eq.CCKP.3.7}
\PP\left[ {\bf D}(T) \in x\,A  \right]& =\PP\left[ \sum_{i=1}^{N(T)} {\bf X}^{(i)}\, e^{-\xi(\tau_i)} \in x\,A \right]\nonumber\\[2mm]
 &= \left(\sum_{n=1}^{M} +  \sum_{n=M+1}^{\infty}\right)\PP\left[ \sum_{i=1}^{n} {\bf X}^{(i)}\, e^{-\xi(\tau_i)} \in x\,A\,,\;N(T)=n \right]\nonumber\\[2mm]
&=:J_1(T,\,x,\,M) + J_2(T,\,x,\,M)\,.
\end{align}

Let start with estimation of $J_2(T,\,x,\,M)$. From \cite[Prop. 2.4]{konstantinides:passalidis:2024g}, taking into account that the process $\{\xi(t)\,,\;t\geq 0\}$ has non-decreasing sample paths and according to Lemma \ref{lem.CCKP.3.2}, for any $\varepsilon>0$, there exists some positive constant $C$ such that
\begin{align}\label{CZT.4.7.0722}
&J_2(T,\,x,\,M) \leq \sum_{n=M+1}^{\infty} \PP\left[ \sum_{i=1}^{n} Y_A^{(i)}\, e^{-\xi(\tau_1)}> x\,,\;N(T)=n \right]\nonumber \\[2mm] 
&\leq C\,\sum_{n=M+1}^{\infty} (1+\vep)^n\, \sum_{i=1}^{n} \PP\left[ Y_A^{(i)}\, e^{-\xi(\tau_1)}> x,\, N(T)=n\right]\nonumber   \\[2mm] 
&\leq C\,\sum_{n=M+1}^{\infty} (1+\vep)^n\,n\,\int_0^T \PP\left[ {\bf X}\, e^{-\xi(s)}\in  x\,A \right]\, \PP[N(T-s)=n-1]\,\PP[\tau_1 \in ds]\nonumber\\[2mm] 
&= C\,\int_0^T\,\sum_{n=M+1}^{\infty} (1+\vep)^n\,n\,  \PP[N(T-s)=n-1]\,\PP\left[ {\bf X}\, e^{-\xi(s)}\in  x\,A \right]\,\PP[\tau_1 \in ds] \nonumber \\[2mm] 
& =C\,\int_0^T \E\left[ (1+\vep)^{N(T-s)+1}\,(N(T-s)+1)\,{\bf 1}_{\{N(T-s)\geq M\}} \right]\,\PP\left[ {\bf X}\, e^{-\xi(s)}\in  x\,A \right]\,\PP[\tau_1 \in ds]\nonumber\\[2mm] 
&\leq C\, \E\left[ (1+\vep)^{N(T)+1}\,(N(T)+1)\,{\bf 1}_{\{N(T)\geq M\}} \right]\,\int_0^T \PP\left[ {\bf X}\, e^{-\xi(s)}\in  x\,A \right]\,\lambda(ds)\,,
\end{align}
where we used the Fubini's theorem in the fourth step.
Due to the fact that $\{N(t)\,,\;t\geq 0\}$ is a renewal process, by \cite{stein:1946} or \cite[Th. 1]{kocetova:leipus:siaulys:2009}, for any $T\in\varLambda\setminus\{\infty\}$, there exists some $\vep_T>0$, such that 
\beao
\E\left[(1+\vep)^{N(T)} \right] < \infty\,
\eeao 
holds for all $0<\vep<\vep_T$. Hence,
\beam \label{eq.CCKP.3.8}
\lim_{M \to \infty} \limsup \dfrac{J_2(T,\,x,\,M)}{\int_0^T \PP\left[  {\bf X}\,e^{-\xi(s)} \in x\,A \right]\,\lambda(ds) } =0\,.
\eeam

Next, we handle the $J_1(T,\,x,\,M)$. Now, we find
\beam \label{eq.CCKP.3.9} \notag
&&J_1(T,\,x,\,M) =\sum_{n=1}^{M}\,\PP\left[ \sum_{i=1}^{n} {\bf X}^{(i)}\, e^{-\xi(\tau_i)}{\bf 1}_{(N(T)=n)} \in x\,A\right]\\[2mm] \notag
&&\sim \sum_{n=1}^{M} \sum_{i=1}^{n}\PP\left[ {\bf X}^{(i)}\, e^{-\xi( \tau_i )} \in x\,A,\, N(T)=n\right]\\[2mm] \notag
&&=\left(\sum_{n=1}^{\infty} -\sum_{n=M+1}^{\infty} \right)\,\sum_{i=1}^{n}\PP\left[ {\bf X}^{(i)}\, e^{-\xi(\tau_i)} \in x\,A\,,\;N(T)=n\right]\\[2mm]
&&=: J_{11}(T,\,x,\,M)-J_{12}(T,\,x,\,M)\,,
\eeam
where in the second step we used Lemma \ref{lem.CCKP.3.1}. For the first term, after changing the order of summation, it holds that
\begin{align} \label{eq.CCKP.3.10} \notag
J_{11}(T,\,x,\,M) &=\sum_{n=1}^{\infty}\,\PP\left[ {\bf X}^{(i)}\, e^{-\xi(\tau_i)} \in x\,A\,,\;\tau_i \leq T \right]\\[2mm]
&=\int_0^T\PP\left[ {\bf X}\, e^{-\xi(s)}\in  x\,A \right]\,\lambda(ds)\,.
\end{align}
For the second term $J_{12}(T,\,x,\,M)$, by the same technique used in the derivation of \eqref{CZT.4.7.0722}, we obtain
\beao
J_{12}(T,\,x,\,M) &\leq& \sum_{n=M+1}^{\infty}\,\sum_{i=1}^{n}\PP\left[Y_A^{(i)}\, e^{-\xi(\tau_1)} > x\, ,\,N(T)=n\right]\\[2mm] \notag
&\leq& \E\left[ [N(T)+1]\,{\bf 1}_{\{N(T) \geq M\}}\right] \int_0^T \PP\left[ {\bf X}\, e^{-\xi(s)} \in x\,A\right]\,\lambda( ds)\,.
\eeao
Since the counting process $\{N(t)\,,\;t\geq 0\}$ is a renewal process with a finite renewal function, therefore, it holds that
\beam \label{eq.CCKP.3.11} 
\lim_{M \to \infty} \limsup \dfrac{J_{12}(T,\,x,\,M)}{\int_0^T \PP\left[  {\bf X}\,e^{-\xi(s)} \in x\,A \right]\,\lambda(ds) } =0\,.
\eeam
By relations \eqref{eq.CCKP.3.10} and \eqref{eq.CCKP.3.11} in combination with relation \eqref{eq.CCKP.3.9}, we obtain that it holds that
\beam \label{eq.CCKP.3.12} 
J_{1}(T,\,x,\,M) \sim \int_0^T \PP\left[ {\bf X}\, e^{-\xi(s)}\in  x\,A \right]\,\lambda( ds)\,.
\eeam
Hence, putting in relation \eqref{eq.CCKP.3.7} the results from \eqref{eq.CCKP.3.8} and \eqref{eq.CCKP.3.12}, we obtain that \eqref{eq.CCKP.1.3} is true.
~\halmos



\subsection{Proof of Theorem \ref{th.CCKP.4.1}}\label{sec.CCKP.4.2} 
Likewise, we need a preliminary lemma that is related to the ``insensitivity" of the multivariate linear single big jump principle of the scaled mixture sums with respect to  ${QAI}_A$, under some finite moments conditions for random weights. This lemma extends \cite[Th. 4.1(i)]{konstantinides:passalidis:2024g} in weighted case.

\ble \label{lem.CCKP.4.1}
Let $A \in \mathscr{R}$ be a fixed set. We suppose that the  ${\bf Z}^{(1)},\,\ldots,\,{\bf Z}^{(n)}$ are  ${QAI}_A$, with distributions $V_1,\,\ldots,\,V_n \in \mathcal{C}_A$, respectively. Let $W_1,\,\ldots,\,W_n$ be $n$ arbitrarily dependent nonnegative random variables (nondegenerate to zero), which are independent of ${\bf Z}^{(1)},\,\ldots,\,{\bf Z}^{(n)}$. If there exist some $p> \bigvee_{i=1}^n J_{V_A^{(i)}}^+$, where $V_A^{(i)}$ is given by \eqref{eq.CCKP.2.15}, such that for any $i=1,\,\ldots,\,n$,
\beam \label{eq.CCKP.4.0}
\E\left[ W_i^p \right]< \infty\,,
\eeam
holds, then relation \eqref{eq.CCKP.3.3} is true.
\ele

\noindent\pr~
The upper bound of \eqref{eq.CCKP.3.3} follows, as with relation (\ref{eq.CZT.3.1}), by applying \cite[Th. 3.2]{chen:yuen:2009}. The lower bound of \eqref{eq.CCKP.3.3} follows, as with relation  \eqref{eq.CCKP.3.5}, where in the last step we use, for any  $i=1,\,\ldots,\,n$, the following relation
\begin{align} \label{eq.CCKP.4.1} 
&\quad \PP\left[ W_i\,{\bf Z}^{(i)} \in  x\,A\,,\; W_j\,{\bf Z}^{(j)} \in  x\,A\right]= \PP\left[ W_i\, Z_A^{(i)} > x\,,\; W_j\,Z_A^{(j)} >x \right] \\[2mm] \nonumber
& =o\left(\PP\left[ W_i\, Z_A^{(i)} > x \right]+\PP\left[ W_j\,Z_A^{(j)} >x \right] \right)=o\left(\PP\left[ W_i\,{\bf Z}^{(i)} \in  x\,A\right]+\PP\left[ W_j\,{\bf Z}^{(j)} \in  x\,A\right] \right),
\end{align}
where in the second step of \eqref{eq.CCKP.4.1}, we used \cite[Lem. 3.1]{chen:yuen:2009}, see also in \cite[Th. 2.2]{li:2013}, due to the fact that $V_A^{(i)} \in \mathcal{C}$ and condition \eqref{eq.CCKP.4.0}.
~\halmos

Now we proceed to the proof of Theorem \ref{th.CCKP.4.1}.

\noindent{\bf Proof of Theorem \ref{th.CCKP.4.1}}~
Let $M$ be some large integer. Then, by the total probability theorem, we obtain that
\beam \label{eq.CCKP.4.2} 
&&\PP\left[ {\bf D}(T) \in  x\,A\right] =\left(\sum_{n=1}^M + \sum_{n=M+1}^{\infty} \right) \PP\left[ \sum_{i=1}^n {\bf X}^{(i)}\,e^{-\xi(\tau_i)} \in x\,A\,,\;N(T)=n \right] \\[2mm] \notag
&&=: I_1(T,\,x,\,M) +  I_2(T,\,x,\,M)\,.
\eeam
Let us start with the second term $I_2(T,\,x,\,M)$. For some $p> J_{F_A}^+\,$, and large enough $x_0>0$ we have
\beam \label{eq.CCKP.4.3}  \notag
&& I_2(T,\,x,\,M) \leq \sum_{n=M+1}^{\infty} \PP\left[ \sum_{i=1}^n Y_A^{(i)}\,e^{-\xi(\tau_i)} > x\,,\;N(T)=n \right] \\[2mm] \notag
&&=\left( \sum_{n=M+1}^{\left\lfloor x/x_0\right\rfloor}+\sum_{n=\left\lfloor x/x_0\right\rfloor+1}^{\infty}\right) \PP\left[ \sum_{i=1}^n Y_A^{(i)}\,e^{-\xi(\tau_i)} > x\,,\;N(T)=n \right] \\[2mm] \notag
&& \leq \sum_{n=M+1}^{\left\lfloor x/x_0\right\rfloor} \PP\left[ \sum_{i=1}^n Y_A^{(i)}e^{K_T}> x \right]\,\PP[N(T)=n] + \PP[N(T)> \left\lfloor x/x_0\right\rfloor]\\[2mm] \notag
&& \leq \sum_{n=M+1}^{\left\lfloor x/x_0\right\rfloor} n\,\PP\left[ Y_A e^{K_T} > \dfrac xn \right]\,\PP[N(T)=n] + \left\lfloor x/x_0\right\rfloor^{-(p+1)}\,\E\left[ N^{p+1}(T)\,{\bf 1}_{\{N(T) > \left\lfloor x/x_0\right\rfloor\}} \right]\\[2mm] \notag
&& =O\left(\PP\left[ Y_A > x \right]\,\E\left[N^{p+1}(T)\,{\bf 1}_{\{M< N(T) \leq x\}} \right] + \left\lfloor x/x_0\right\rfloor^{-(p+1)}\,\E\left[ N^{p+1}(T)\,{\bf 1}_{\{N(T) > \left\lfloor x/x_0\right\rfloor\}} \right]\right)\\[2mm]
&&=O(\PP\left[ {\bf X} \in x\,A \right]\,\E\left[N^{p+1}(T)\,{\bf 1}_{\{ N(T) > M\}} \right])\,,
\eeam
where in the first step we use \cite[Prop. 2.4]{konstantinides:passalidis:2024g}, in the fourth step we employed Markov's inequality, in the fifth step we used \cite[Lem. 1]{li:2018}, and in the last step we used relation \eqref{eq.CCKP.2.8}. Further, by  \cite[Th. 3.3 (iv)]{cline:samorodnitsky:1994}, Assumption \ref{ass.CCKP.1.2}, and Fatou's lemma, we obtain that
\beam \label{eq.CCKP.4.4}  \notag
&&\liminf\int_0^T\dfrac { \PP\left[  {\bf X}\,e^{-\xi(s)} \in x\,A \right]}{\PP\left[ {\bf X} \in x\,A \right]}\lambda(ds)\\[2mm]
&& \geq \int_0^T\liminf\dfrac { \PP\left[ {Y_A}\,e^{-\xi(s)} >x \right]}{\PP\left[ Y_A  > x \right]}\lambda(ds) = \int_0^T \E\left[\overline{F_A}_{*}\left(e^{\xi(s)}\right)\right]\lambda(ds)>0\,.
\eeam
Thus, by relations \eqref{eq.CCKP.4.3} and \eqref{eq.CCKP.4.4} we find that 
\beam \label{eq.CCKP.4.5} 
&& \lim_{M \to \infty} \limsup \dfrac{I_2(T,\,x,\,M)}{\int_0^T \PP\left[  {\bf X}\,e^{-\xi(s)} \in x\,A \right]\,\lambda(ds) } =0\,.
\eeam

Let us proceed now to the estimate of $I_1(T,\,x,\,M)$. 
\begin{align} \label{eq.CCKP.4.6}  \notag
&I_1(T,\,x,\,M) = \sum_{n=1}^M \PP\left[ \sum_{i=1}^n {\bf X}^{(i)}\,e^{-\xi(\tau_i)} \,{\bf 1}_{\{N(T)=n\}}  \in x\,A \right] \\[2mm]
& \sim \sum_{n=1}^M \sum_{i=1}^n \PP\left[ {\bf X}^{(i)}\,e^{-\xi( \tau_i)} \in x\,A,\, N(T)=n\right]\\[2mm] \notag
&=\left( \sum_{n=1}^{\infty} -  \sum_{n=M+1}^{\infty} \right)\,\sum_{i=1}^n \PP\left[ {\bf X}^{(i)}\,e^{-\xi(\tau_i)} \in x\,A\,,\; N(T)=n \right] =: I_{11}(T,\,x,\,M) - I_{12}(T,\,x,\,M)\,,
\end{align}
where in the second step we used Lemma \ref{lem.CCKP.4.1}. 

Further, for the first term $I_{11}(T,\,x,\,M)$ we change the order of summation to find the expression
\beam \label{eq.CCKP.4.7} 
I_{11}(T,\,x,\,M)= \sum_{i=1}^\infty \PP\left[ {\bf X}\,e^{-\xi(\tau_i)} \in x\,A\,,\; \tau_i \leq T \right]=\int_0^T \PP\left[ {\bf X}\,e^{-\xi(s)} \in x\,A\right] \,\lambda(ds)\,.
\eeam
For the other term $I_{12}(T,\,x,\,M)$, we find the inequalities
\beao
I_{12}(T,\,x,\,M) &\leq& \sum_{n=M+1}^{\infty} \sum_{i=1}^n \PP\left[{\bf X}^{(i)} e^{K_T} \in x\,A\right]\,\PP[N(t)=n]\\[2mm] 
&=& O\left(\PP[{\bf X} \in x\,A] \E\left[ N(T)\,{\bf 1}_{\{ N(T) > M\}}\right]\right)\,.
\eeao 
Hence, by relation \eqref{eq.CCKP.4.4} we obtain that
\beam \label{eq.CCKP.4.8} 
&&\lim_{M \to \infty} \limsup \dfrac{I_{12}(T,\,x,\,M)}{\int_0^T \PP\left[ {\bf X}\,e^{-\xi(s)} \in x\,A\right] \,\lambda(ds)}=0.
\eeam
Therefore, from relations \eqref{eq.CCKP.4.7} and \eqref{eq.CCKP.4.8}, in combination with relation \eqref{eq.CCKP.4.6}, it holds that
\beam \label{eq.CCKP.4.9} 
I_{1}(T,\,x,\,M) \sim \int_0^T \PP\left[ {\bf X}\,e^{-\xi(s)} \in x\,A\right] \,\lambda(ds)\,.
\eeam
As a result, by relations \eqref{eq.CCKP.4.5} and \eqref{eq.CCKP.4.9}, in combination with relation  \eqref{eq.CCKP.4.2} we find that  \eqref{eq.CCKP.1.3} holds. 
~\halmos

Next, we prove Corollary \ref{cor.CCKP.4.1}. 

\noindent{\bf Proof of Corollary \ref{cor.CCKP.4.1}}~
It is sufficient to prove 
\begin{align}\label{coro4.22}
	\int_0^T \PP\left[ {\bf X}\,e^{-\xi(s)} \in x\,A\right] \,\lambda(ds)\sim \mu(A)\overline{G}(x)\int_0^T \E\left[ e^{-\alpha\,\xi(s)}\right]\,\lambda(ds).
\end{align}
First, by Assumption \ref{ass.CCKP.1.2} we can apply to Breiman's theorem. Hence for any $s\in [0,T]$, it holds that 
\begin{align*}
	\PP\left[ {\bf X}\,e^{-\xi(s)} \in x\,A\right]&=\PP\left[ {Y}_A\,e^{-\xi(s)} >x\right]\sim \E\left[e^{-\alpha\xi(s)}\right]\overline{F_A}(x)\\[2mm]
	&=\E\left[e^{-\alpha\xi(s)}\right]\PP\left[ {\bf X} \in x\,A\right] \sim\mu(A) \E\left[e^{-\alpha\xi(s)}\right]\overline{G}(x).
\end{align*}
Nevertheless, by Assumption \ref{ass.CCKP.1.2} and \cite[Prop. 2.2.1]{bingham:goldie:teugels:1987}, we have 
\begin{align*}
	\frac{\PP\left[ {\bf X}\,e^{-\xi(s)} \in x\,A\right]}{\PP\left[ {\bf X} \in x\,A\right]}\leq\frac{\overline{F_A}\left(e^{-K_T}x\right)}{\overline{F_A}(x)}\leq 2e^{-2\alpha K_T}.
\end{align*} 
Obviously, $2e^{-2\alpha K_T}$ is integrable, then \eqref{coro4.22} is verified by the dominated convergence theorem.
~\halmos

\subsection{Proof of Theorem \ref{th.CCKP.5.1} }\label{subsec.CCKP.5.2}

Before proving  Theorem \ref{th.CCKP.5.1}, we need a preliminary lemma. This lemma follows after some elementary modifications in the proof of \cite[Th. 2]{yi:chen:su:2011}. However, for the sake of completeness, we provide a full proof with a slightly different way.

\ble \label{lem.CCKP.5.1}
Under the conditions of Theorem \ref{th.CCKP.5.1}, it holds
\beam \label{eq.CCKP.5.6} 
\lim_{M \to \infty}\limsup \dfrac{\sum_{i=M+1}^{\infty}\PP\left[ {\bf X}^{(i)}\ e^{-\xi(\tau_i)} \in x\,A \right]}{\PP\left[ {\bf X}\ e^{-\xi(\tau_1)} \in x\,A \right]}=0\,.
\eeam  
\ele

\noindent\pr~
Let  $M,\,x$ be large enough, then due to \cite[Lem. 1]{yi:chen:su:2011}, then for $0<p_1 < J_{F_A}^- \leq J_{F_A}^+ < p_2 <\infty $, there exist positive constants $C_1$, $C_2$ (both irrespective of $e^{-\xi(\tau_i)}$) such that
\beam \label{eq.CCKP.5.7} \notag
&&\PP\left[ {\bf X}^{(i)} e^{-\xi(\tau_i)} \in x\,A \right] \leq C_1\,\PP\left[ {\bf X} \in x\,A \right]\,\left( \E\left[e^{-p_1\,\xi(\tau_i)} \right] \bigvee \E\left[e^{-p_2\,\xi(\tau_i)} \right] \right)\\[2mm]
&&\leq C_2\, \PP\left[ {\bf X}\ e^{-\xi(\tau_1)}\in x\,A \right]\,\left( \E\left[e^{-p_1\,\xi(\tau_i)} \right] \bigvee \E\left[e^{-p_2\,\xi(\tau_i)} \right] \right)^{\frac1{\rho}}\,,
\eeam
where in the second step we employed \cite[Th. 3.3(iv)]{cline:samorodnitsky:1994}.
Through relation \eqref{eq.CCKP.5.7}, it holds that
\begin{align*}
&\lim_{M\to\infty}\limsup \dfrac{\sum_{i=M+1}^{\infty}\PP\left[ {\bf X}^{(i)} e^{-\xi(\tau_i)} \in x\,A \right]}{\PP\left[ {\bf X}\ e^{-\xi(\tau_1)} \in x\,A \right]}\\[2mm]
& \leq C_2\, \lim_{M\to\infty}\left[ \sum_{i=M+1}^{\infty} \,\left( \E\left[e^{-p_1\,\xi(\tau_i)} \right] \bigvee \E\left[e^{-p_2\,\xi(\tau_i)} \right] \right)^{\frac1{\rho}}\right]=0 \,.
\end{align*}
This concludes \eqref{eq.CCKP.5.6}.

~\halmos

\noindent{\bf Proof of Theorem \ref{th.CCKP.5.1}}~
For the upper bound, we obtain that
\beam \label{eq.CCKP.5.9} \notag
&&\PP[{\bf D}(\infty) \in x\,A] =\PP\left[ \sum_{i=1}^{\infty} {\bf X}^{(i)} e^{-\xi(\tau_i)} \in x\,A \right]=\PP\left[ \sup_{{\bf p} \in I_A} {\bf p}^T\left(\sum_{i=1}^{\infty} {\bf X}^{(i)} e^{-\xi(\tau_i)} \right) > x \right]\\[2mm] \notag
&&\leq\PP\left[ \sum_{i=1}^{\infty} Y_A^{(i)} e^{-\xi(\tau_i)}> x \right]\sim \sum_{i=1}^{\infty} \PP\left[ Y_A^{(i)} e^{-\xi(\tau_i)}> x \right]= \sum_{i=1}^{\infty} \PP\left[{\bf X} e^{-\xi(\tau_i)}\in  x\,A\,\right]\\[2mm]
&&= \int_{0}^{\infty} \PP\left[ {\bf X}\,e^{-\xi(s)}  \in x\,A \right]\,\lambda(ds)\,.
\eeam
where in the fourth step, we used \cite[Th. 2]{yi:chen:su:2011}. 

Next, we deal with the lower bound. Let $M \in \bbn$, be large enough, then it holds
\beam \label{eq.CCKP.5.10}  \notag
&&\PP[{\bf D}(\infty) \in x\,A] \geq \PP\left[ \sum_{i=1}^{M} {\bf X}^{(i)} e^{-\xi(\tau_i)} \in x\,A \right]  \geq  \PP\left[ \bigcup_{i=1}^{M}\left\{ {\bf X}^{(i)} e^{-\xi(\tau_i)} \in x\,A \right\} \right]\\[2mm] \notag
&&\geq \sum_{i=1}^{M} \PP\left[ {\bf X}^{(i)} e^{-\xi(\tau_i)} \in x\,A \right] - \sum_{1\leq i < j \leq M}  \PP\left[ {\bf X}^{(i)} e^{-\xi(\tau_i)} \in x\,A\,,\; {\bf X}^{(j)} e^{-\xi(\tau_j)} \in x\,A\right] \\[2mm] \notag
&&\sim \sum_{i=1}^{M} \PP\left[ {\bf X}^{(i)} e^{-\xi(\tau_i)} \in x\,A \right]=\left(\sum_{i=1}^{\infty} - \sum_{i=M+1}^{\infty} \right)\PP\left[ {\bf X}^{(i)} e^{-\xi(\tau_i)} \in x\,A \right]\\[2mm]
&&\sim  \int_{0}^{\infty} \PP\left[ {\bf X}\,e^{-\xi(s)}  \in x\,A \right]\,\lambda(ds)\,.
\eeam
where at the second step we used the fact that $x\,A$ is increasing, at the third step we used the Bonferroni inequality, and at the fourth step we used \cite[Th. 2.2]{li:2013}, while in the last step we let $M \to \infty$ and we applied Lemma \ref{lem.CCKP.5.1}. From relations \eqref{eq.CCKP.5.9} and \eqref{eq.CCKP.5.10} we obtain the desired result.
~\halmos

\noindent{\bf Proof of Corollary \ref{cor.CCKP.5.1}}~
Analogously to the proof of Corollary \ref{cor.CCKP.4.1}, we only need to prove \eqref{coro4.22} with $T=\infty$. The convergence is obvious; we only need to find the dominant function. Due to \cite[Lem. 1]{yi:chen:su:2011}, for any $s\geq 0$ and some fixed $0<p_1 < \alpha< p_2 <\infty$ indicated in Assumption \ref{ass.CCKP.5.1}, there exists some constant $C>0$ such that
\begin{align*}
	\frac{\PP\left[ {\bf X}\,e^{-\xi(s)}  \in x\,A \right]}{\PP\left[ {\bf X}  \in x\,A \right]}&\leq C\E\left[e^{-p_1\xi(s)}\right]\vee \E\left[e^{-p_2\xi(s)}\right],
\end{align*}
which is integrable because of \eqref{eq.CCKP.5.2}. Thus, Corollary \ref{cor.CCKP.5.1} is proved by the dominated convergence theorem.
~\halmos

\section*{Acknowledge}
The authors are grateful to Prof. Dongya Cheng for her meticulous review of this manuscript and for providing Example \ref{eg2.2}, which significantly enhanced the clarity of our work.

\noindent \textbf{Disclosure statement.}
There are no financial or non-financial competing interests.

\end{document}